%% file: main_arxiv.tex
\documentclass{article}
\usepackage{authblk}  

\usepackage{PRIMEarxiv}

\usepackage{graphicx} 
\usepackage[table]{xcolor}
\usepackage{amsfonts}
\usepackage{amsmath}
\usepackage{amssymb}
\usepackage{amstext}
\usepackage{amsthm}
\usepackage{bm}
\usepackage{url}
\usepackage{subcaption}
\usepackage{stmaryrd}
\usepackage{multirow}
\usepackage{esint} 
\usepackage{comment}
\usepackage[colorlinks=true,urlcolor=blue,linkcolor=blue,citecolor=blue]{hyperref}
\usepackage{booktabs}
\usepackage{makecell}
\usepackage{mathtools}

\usepackage[capitalise,noabbrev]{cleveref}
\crefformat{equation}{(#2#1#3)}
\numberwithin{equation}{section}
\numberwithin{figure}{section}
\usepackage[normalem]{ulem}
\newtheorem{remark}{Remark}[section]

\newcommand{\tA}{\tilde{A}}
\newcommand{\tB}{\tilde{B}}
\newcommand{\Ntr}{N_{\rm tr}}
\newcommand{\tspan}{\text{span}}

\DeclareMathOperator*{\argmax}{arg\,max}

\input ./pf

\title{Symbolic Regression of Data-Driven Reduced Order Model Closures for Under-Resolved, Convection-Dominated Flows}

\author[1]{Simone Manti}
\author[2]{Ping-Hsuan Tsai}
\author[1]{Alessandro Lucantonio}
\author[2]{Traian Iliescu}

\affil[1]{Aarhus University, Katrinebjergvej 89G-F, 8000, Aarhus, Denmark}
\affil[2]{Department of Mathematics, Virginia Tech, 225 Stanger Street, Blacksburg, 24061, VA, USA}

\begin{document}

\footnotetext{Corresponding author: \texttt{a.lucantonio@mpe.au.dk} (Alessandro Lucantonio)}

\maketitle

\begin{abstract}
   Data-driven closures correct the standard reduced order models (ROMs) to increase their accuracy in under-resolved, convection-dominated flows. There are two types of data-driven ROM closures in current use: (i) structural, with simple ansatzes (e.g., linear or quadratic); and (ii) machine learning-based, with neural network ansatzes. We propose a novel symbolic regression (SR) data-driven ROM closure strategy, which combines the advantages of current approaches and eliminates their drawbacks. As a result, the new data-driven SR closures yield ROMs that are interpretable, parsimonious, accurate, generalizable, and robust. To compare the data-driven SR-ROM closures with the structural and machine learning-based ROM closures, we consider the data-driven variational multiscale ROM framework and two under-resolved, convection-dominated test problems: the flow past a cylinder and the lid-driven cavity flow at Reynolds numbers $\rm{Re} = 10000, 15000, \text{and } 20000$. This numerical investigation shows that the new data-driven SR-ROM closures yield more accurate and robust ROMs than the structural and machine learning ROM closures.
\end{abstract}


\section{Introduction}

High-performance computing and modern numerical algorithms have made high-fidelity fluid-thermal analysis tractable in geometries of ever increasing complexity. Despite continued advances in these areas, direct numerical (DNS), large eddy simulation (LES), and even unsteady Reynolds-averaged Navier–Stokes (URANS) simulations of turbulent thermal transport remain too costly for routine analysis and design of thermal–hydraulic systems, where hundreds of cases must be considered. 

Reduced order models (ROMs) offer a promising alternative by leveraging expensive
high-fidelity simulations (referred to as full order models or FOMs) to first extract a low-dimensional basis that captures the principal features of the underlying flow fields, {and then construct computational models whose dimensions are orders of magnitude lower than the FOM dimension}.
In the numerical simulation of fluid flows, Galerkin ROMs (G-ROMs), which use data-driven basis functions in a Galerkin framework, have provided efficient and accurate approximations of laminar flows, such as the two-dimensional flow past a circular cylinder at low Reynolds
numbers~\cite{hesthaven2015certified,
quarteroni2015reduced}. However, turbulent flows are notoriously hard for the standard G-ROM. 
Indeed, to capture the complex dynamics, a large number \cite{tsai2023accelerating} of ROM basis functions is required, which yields high-dimensional ROMs that cannot be used in realistic applications.
Thus, computationally efficient, low-dimensional ROMs are used instead.
Unfortunately, these ROMs are inaccurate since the ROM basis functions that were not used to build the G-ROM have an important role in dissipating the energy from the system~\cite{ahmed2021closures}.
Indeed, without enough dissipation, the low-dimensional G-ROM generally yields spurious numerical oscillations.
Thus, closures and stabilization strategies are required for the low-dimensional G-ROMs to be stable and accurate~\cite{ahmed2021closures,
fick2018stabilized,kaneko2022augmented,kaneko2020towards,parish2024residual,tsai2024time}. 

The ROM closures are ROM corrections in which the standard G-ROM is supplemented with a ROM closure model, $\boldsymbol{\tau}^{\text{ROM}}$, which models the effect of the unresolved scales (i.e., the effect of the ROM basis functions that were not used to build the G-ROM).  
There are several types of ROM closure modeling strategies, which are carefully reviewed in \cite{ahmed2021closures}.

Arguably the most popular class of ROM closures is data-driven ROM closures
\cite{ahmed2021closures,sanderse2024scientific}, in which FOM data
are used to train the ROM closure model. There are several 
data-driven approaches for modeling $\boldsymbol{\tau}^{\text{ROM}}$, which are surveyed in~\cite{ahmed2021closures}.
To construct the data-driven ROM closure model, $\boldsymbol{\tau}^{\text{ROM}}$, we first 
postulate a model form for $\boldsymbol{\tau}^{\text{ROM}}$. 
Currently, there are two strategies for the ROM closure model form: %
\begin{enumerate}
    \item A principled, structural strategy, in which the ROM closure model ansatz has a simple mathematical structure (e.g., 
    a linear model or a quadratic model~\cite{mou2021data,xie2018data}).
    \item A nonintrusive, machine learning (ML) strategy, in which the ROM closure model ansatz is based on neural networks~\cite{ahmed2021closures,sanderse2024scientific,ahmed2023physics,xie2020closure}.
\end{enumerate}
In the final step of the construction of the data-driven ROM closure model, we leverage the available data to determine the unknown parameters in the ROM closure model ansatz.
This can be achieved, e.g.,  by solving a least squares problem (for the structural strategy) or by minimizing a loss function (for the machine learning strategy).

Each strategy has pros and cons:
\begin{enumerate}
    \item The structural strategy yields simple, interpretable ROM closure models, which can be easily combined with the ROM operators in the standard G-ROM.
    However, it is not clear that the simple building blocks (e.g., the linear and quadratic terms) can accurately represent the complex dynamics modeled by the ROM closure terms.
    \item The machine learning strategy, on the other hand, yields models that can, in principle, approximate the complex underlying dynamics. However, this strategy yields ROM closure models that contain parameters whose meaning is generally unclear.
\end{enumerate}

Finding data-driven ROM closure models that are both accurate and physically interpretable is a major challenge in the development of ROMs for realistic turbulent flows.
In this paper, we propose a novel data-driven ROM closure modeling strategy that centers around symbolic regression (SR).
The novel SR-ROM closure modeling strategy combines the advantages of the current approaches and eliminates their drawbacks: The SR-ROM closure models are interpretable, parsimonious (i.e., contain 
much fewer parameters than, e.g., the machine learning strategy), and are rich enough to approximate the complex ROM closure dynamics.
Given the limited number of parameters that characterize the new SR-ROM closures, they are less prone to overfitting and, thus, generalize better than the current data-driven ROM closures.
Furthermore, the new SR-ROM closure models have the following additional advantages:
They are accurate, computationally efficient, and, because they are constructed by using a sound mathematical strategy, they are amenable to analysis. We emphasize that, although SR has been successfully used in data-driven modeling, to our knowledge, the SR-ROM closure model proposed in this paper is the first instance in which SR is used in ROM closure modeling. 

To illustrate our new SR-ROM closure modeling strategy, we need to choose a data-driven ROM closure modeling approach.
To this end, we choose the data-driven variational multiscale ROM (d2-VMS-ROM) \cite{mou2021data,xie2018data}, which is a successful data-driven ROM closure modeling approach utilized in various applications (from flow past a cylinder to flow past a backward-facing step to the quasi-geostrophic equations\cite{mou2021data,mou2021reduced}). 
The d2-VMS-ROM is a data-driven ROM closure that leverages the intrinsic hierarchical nature of the ROM basis and the VMS framework \cite{hughes1995multiscale} to construct data-driven ROM closures.
Specifically, the d2-VMS-ROM strategy proceeds by decomposing the ROM space $\bX^R \coloneqq \text{span} \{ \bphi_1, \ldots, \bphi_R \}$ into the low-dimensional space of large, resolved scales, $\bX^r \coloneqq \text{span} \{ \bphi_1, \ldots, \bphi_r \}$, and the large-dimensional space of small, unresolved scales, $\bX^{'} \coloneqq \text{span} \{ \bphi_{r+1}, \ldots, \bphi_R \}$.
To ensure computational efficiency, the d2-VMS-ROM is constructed in the low-dimensional space $\bX^r$.
However, to ensure numerical accuracy, the d2-VMS-ROM is supplemented with a ROM closure model, $\boldsymbol{\tau}^{\text{ROM}}$, which models the effect of the unresolved scales in $\bX^{'}$.
We note that one distinguishing feature of the d2-VMS-ROM is that it is built by displaying the
explicit formula for the ROM closure term $\boldsymbol{\tau}^{\text{ROM}}$, which is achieved by leveraging the VMS formulation.  
Once we determine the precise form of the ROM closure term, we use one of the two data-driven modeling strategies outlined above (i.e., either a structural or an ML approach) to determine the unknown parameters of the closure model.

In this paper, we consider several d2-VMS-ROM formulations, which are based on linear regression, quadratic manifold formulation, and ML.
Furthermore, we also consider a d2-VMS-ROM formulation based on the novel SR-ROM closure model.
To assess the performance of the new SR-ROM closure, we perform an investigation of all four formulations in the numerical simulation of two convection-dominated test problems: flow past a cylinder and lid-driven cavity flow.

The rest of the paper is organized as follows: We begin by introducing the FOM in Section~\ref{sec:fom}. The ROM methodology is presented in Section~\ref{subsection:g-rom}. Section~\ref{subsection:d2-vms-roms} provides the relevant background information and notation for the VMS approach. We also introduce
the novel SR-ROM method along with  various ML techniques used to model the closure term for the VMS-ROM, which we refer to as ML-VMS-ROMs. In Section~\ref{sec:nr}, the 
SR-ROM is compared with the current ML-VMS-ROMs in the numerical simulation of two test problems in the under-resolved, convection-dominated regime: the 2D flow past a cylinder and the 2D lid-driven cavity at Reynolds numbers $\rm{Re} = 10000, 15000, \text{and } 20000$. Concluding remarks are presented in Section~\ref{sec:conclusions}.

\section{Full Order Model (FOM)}
\label{sec:fom}

The governing equations are the incompressible Navier-Stokes equations (NSE):
\begin{equation} \label{eq:nse}
  \frac{\partial \bu}{\partial t} + (\bu \cdot \nabla)\bu 
    = -\nabla p + \frac{1}{\rm Re} 
    \Delta \bu, \qquad
   \nabla \cdot \bu = 0.
\end{equation}
Here, $\bu$ represents the velocity subject to {homogeneous} Dirichlet 
boundary conditions, $p$ denotes the pressure,
and $\rm Re$ is the Reynolds number. 

We employ the spectral element method 
for the spatial discretization of
(\ref{eq:nse}). The $P_{q}$--$P_{q-2}$ velocity-pressure coupling
\cite{maday1987well} is considered, where the velocity $\bu$ is represented as a
tensor-product Lagrange polynomial of degree $q$ in the reference element $\Oh
:= [-1,1]^2$, while the pressure $p$ is of degree $q-2$. 
{This yields} a velocity
space $\bX^{\cN}$ for approximating the velocity, and a pressure space
$Y^{\bar{\cN}}$ for approximating the pressure, where the finite dimensions
$\cN$ and $\bar{\cN}$ are the global numbers of spectral element degrees of
freedom in the corresponding spectral element spaces.

The 
{semi-discretization} of
(\ref{eq:nse}) reads: 
{\em Find $(\bu
,~p) \in
(\bX^{\cN}, Y^{\bar{\cN}})$ such that, for all $(\bv,~q) \in 
(\bX^{\cN},Y^{\bar{\cN}})$,}
\begin{align} \label{eq:var1}
   \left( \pp{\bu}{t},\bv \right)  + \biggl( (\bu \cdot \nabla) \bu, \bv \biggr) & = \left( p, \nabla \cdot
   \bv \right) - \frac{1}{\rm Re} \, (\nabla \bu, \nabla \bv),
   \\[1.6ex] \label{eq:var2}
   - \left(\nabla \cdot \bu, q  \right)  &= 0,
\end{align}
where $(\cdot, \cdot)$ denotes 
the $L^2$ inner product. 
The resulting spatially--discretized ordinary differential equations 
are solved using
a semi-implicit BDF$k$/EXT$k$ scheme \cite{fischer2017recent}. We refer to \cite{tsai2022parametric}
for the detailed derivation of the FOM 
and the treatment with the 
nonhomogeneous Dirichlet boundary conditions. 

\section{Reduced Order Models (ROMs)}

We 
present the Galerkin reduced order model (G-ROM) in Section~\ref{subsection:g-rom}. In Section~\ref{subsection:d2-vms-roms}, we first provide background information and establish notation for the VMS-ROM framework. Then, we leverage several ML strategies to construct VMS-ROM closure models. The resulting ROMs, which we collectively call ML-VMS-ROMs, are of four types:
the linear regression VMS-ROM (Section~\ref{subsec: lr-rom});
the data-driven VMS-ROM (Section~\ref{subsec:d2-vms-rom});
the novel symbolic regression VMS-ROM (Section~\ref{subsec:sr-rom}); and
the neural network VMS-ROM (Section~\ref{subsec:nn-rom}).

\subsection{Galerkin ROM (G-ROM)}
\label{subsection:g-rom}

To construct the G-ROM,
we begin by collecting a set of velocity snapshots $\{\bu^k :=
\bu(\bx,t^k)-\bphi_0\}^K_{k=1}$, which correspond to the FOM solutions at
time instances $t^k$, with the subtraction of the zeroth mode, $\bphi_0$.
We then employ the standard proper orthogonal decomposition (POD) procedure
\cite{berkooz1993proper,volkwein2013proper} to construct the reduced basis
functions. The Gramian matrix is formed using the $L^2$ inner product, and
the first $r$ POD basis functions $\{\bphi_i\}^r_{i=1}$ are constructed from the
first $r$ eigenmodes of the Gramian.

The G-ROM is constructed by inserting the ROM basis expansion
\begin{equation} \label{eq:romu}
   \bu_\tr(\bx) = \bphi_0(\bx) + \sum_{j=1}^r u_{\tr,j} \bphi_j(\bx)
\end{equation}
into
(\ref{eq:nse}):
{\em Find $\bu_\tr$
such that, for all $\bv \in \bX_1$,}
\textbf{}\begin{eqnarray}
    && 
    \left(
        \frac{\partial \bu_\tr}{\partial t}, \bv_{i}  
    \right)
    + \frac{1}{\rm Re} \, 
    \left( 
        \nabla \bu_\tr, \nabla \bv_{i}
    \right)
    + \biggl( (\bu_\tr \cdot \nabla) \bu_\tr, \bv_{i}
    \biggr)
    = 0,  
    \label{eq:gromu}
\end{eqnarray}
where $\bX_1 := \text{span} \{\bphi_i\}^r_{i=1}$ is the reduced space.

\begin{remark}
    In this paper, the zeroth mode $\bphi_0$ is set to be the initial condition of the FOM simulation. 
\end{remark}

\begin{remark}
We note that, in the case of fixed geometries, the divergence and
pressure terms drop out of (\ref{eq:gromu}) since the 
ROM basis is
weakly divergence-free. For ROMs that include the pressure approximation, see, e.g.,~\cite{hesthaven2015certified,quarteroni2015reduced,ballarin2015supremizer,decaria2020artificial,noack2005need}. 
\end{remark}

From (\ref{eq:gromu}), a system of differential equations in the coefficients 
of the POD bases, $u_{\tr,j}$, 
is derived:

\begin{align}
    B \frac{d \uu_{\tr}}{dt} & = -\mC (\bar{\uu}_{\tr}) \bar{\uu}_{\tr} -
    \frac{1}{\rm Re} A \bar{\uu}_{\tr}, \label{eq:nse_ode}
\end{align}
where $\uu_\tr \in \mathbb{R}^r$ is the vector 
of POD coefficients $\{ u_{r,j}\}^r_{j=1}$, and           
$\bar{\uu}_{\tr} \in \mathbb{R}^{r+1}$ is the augmented vector that includes the zeroth mode's coefficient. $A$, $B$, and $\mC$ represent the  reduced stiffness, mass, and advection
operators, respectively, with entries
\begin{align} 
    A_{ij} =\int_{\Omega}\nabla \bphi_i :\nabla \bphi_j\, dV,
    \quad B_{ij} =\int_{\Omega}\bphi_i\cdot\bphi_j\, dV,
    \quad \label{eq:Cu} \mC_{ikj}  =\int_{\Omega}\bphi_i\cdot(\bphi_k\cdot\nabla)\bphi_j\,dV.                   
\end{align}

For temporal discretization of (\ref{eq:nse_ode}), a semi-implicit scheme with $k$th-order backward differencing (BDF$k$) and $k$th-order extrapolation (EXT$k$) is considered. The fully discretized reduced system at time $t^l$ is then derived:
\begin{align}
    \left(\frac{\beta_0}{\Delta t} B + \frac{1}{\rm Re} A \right) \,
    \uu_{\tr}^{l+1} & = - \sum^k_{i=1} \alpha_i \left[ \mC ({\uu}^{l-i}_{\tr})
    {\uu}^{l-i}_{\tr} + (C_1 + C_2) \bar{\uu}^{l-i}_{\tr} - \underline{c}_{0} \right] - B
    \sum^k_{i=1}\frac{\beta_i}{\Delta t} \uu^{l-i}_\tr - \frac{1}{\rm Re}
    \underline{a}_0, \label{eq:nse_d1}
\end{align}
where
\begin{align}
    c_{0,i} & = \int_{\Omega}\bphi_i\cdot(\bphi_0\cdot\nabla)\bphi_0\,dV,\quad
    a_{0,i} = \int_{\Omega}\nabla \bphi_i : \nabla \bphi_0 \,dV, \\ 
    C_{1,ij} & = \int_{\Omega}\bphi_i\cdot(\bphi_0\cdot\nabla)\bphi_j\,dV,\quad
    C_{2,ik} = \int_{\Omega}\bphi_i\cdot(\bphi_k\cdot\nabla)\bphi_0\,dV, \label{eq:zeroth_contribution}
\end{align}
for all $i=1,\ldots,r$ and $j,~k=0,\ldots,r$.

\subsection{Machine Learning Variational Multiscale ROMs (ML-VMS-ROMs)} 
\label{subsection:d2-vms-roms}

The variational multiscale (VMS) methods are numerical discretization techniques that greatly enhance the accuracy of classical Galerkin approximations in {under-resolved} simulations. This situation often occurs with coarse meshes or when there are insufficient basis functions available. The VMS framework, originally proposed by Hughes and his colleagues \cite{hughes1998variational}, has significantly influenced various fields in computational mechanics. For comprehensive reviews of its impact, see \cite{ahmed2017review,
codina2018variational}.

In the following, we describe the VMS-ROM framework \cite{mou2021data}. 
{First,} we consider two reduced spaces $\bX_1 \coloneqq
\tspan\{\bphi_1,\ldots,\bphi_r\}$ and $\bX_2 \coloneqq
\tspan\{\bphi_{r+1},\ldots,\bphi_R\}$, where $\bphi_i$ is the $i$-th reduced
basis function, and $R$ is the dimension of the snapshot dataset. We note that
$\bX_1$ represents the span of the resolved ROM scales, and $\bX_2$ represents
the span of the unresolved scales. In addition, the two reduced spaces are
orthogonal because of the POD construction. 

Next, we use the best ROM approximation of $\bu$ in the space $\bX_1 \oplus \bX_2$, that is, $\bu_{\tR} \in \bX_1 \oplus \bX_2$, defined as
\begin{equation}
    \bu_R = 
    \sum^R_{j=1} u_{r,j} \bphi_j(\bx) = \sum^r_{j=1}
    u_{r,j} \bphi_j + \sum^R_{j=r+1} u_{r,j} \bphi_j = \bu_r + \bu'_r, \label{eq:rom_expansion_R}
\end{equation} 
where $\bu_\tr \in \bX_1$ represents
the resolved ROM component of $\bu$, and $\bu' \in \bX_2$ represents the unresolved ROM
component of $\bu$. 

Plugging $\bu_\tR$ in (\ref{eq:nse}) and projecting the resulting equation onto $\bX_1$, we obtain
\begin{eqnarray}
    && 
    \left(
        \frac{\partial \bu_\tR}{\partial t}, \bv_{i}
    \right)
    + \frac{1}{\rm Re} \, 
    \left( 
        \nabla \bu_\tR, \nabla \bv_{i}
    \right)
    + \biggl( 
        (\bu_\tR \cdot \nabla) \bu_\tR, \bv_{i}
    \biggr)
    = 0, \qquad \forall~i=1,\ldots,r. \label{eq:rom_expansion_R_2}
\end{eqnarray}
{Using} the orthogonality of the ROM basis functions, one can show that 
\begin{equation}
    \left( \pp{\bu_\tR}{t}, \bv_i  \right) = \left(  \pp{\bu_\tr}{t}, \bv_i \right), \qquad \forall~i=1,\ldots,r,
\end{equation}
and (\ref{eq:rom_expansion_R_2}) can be further written as: $\forall ~i=1,\ldots,r$
\begin{align}
    \left(
        \frac{\partial \bu_\tr}{\partial t}, \bv_{i}
    \right)
    & + \frac{1}{\rm Re} \, 
    \left( 
        \nabla \bu_\tr, \nabla \bv_{i}
    \right)
    + \biggl( 
        (\bu_\tr \cdot \nabla) \bu_\tr, \bv_{i}
    \biggr) \nonumber \\ & +  
\left[
    \frac{1}{\rm Re} \, 
    \left( 
        \nabla \bu_\tR, \nabla \bv_{i}
    \right)
    + \biggl( 
        (\bu_\tR \cdot \nabla) \bu_\tR, \bv_{i}
    \biggr)  - 
    \frac{1}{\rm Re} \, 
    \left( 
        \nabla \bu_\tr, \nabla \bv_{i}
    \right)
    - \biggl( 
        (\bu_\tr \cdot \nabla) \bu_\tr, \bv_{i}
    \biggr)
    \right] = 0.
    \label{eq:vmsrom}
\end{align}
The bracketed term in \cref{eq:vmsrom} is the VMS-ROM closure term\footnote{If one considers a zeroth mode in the definition of $\bu_R$ and $\bu_r$, the zeroth mode contribution will not appear in the VMS-ROM closure term because of the cancellation. 
}, which models the
interaction between the ROM modes $\{\bphi_1,\ldots,\bphi_r\}$ and the discarded
ROM modes $\{\bphi_{r+1},\ldots,\bphi_R\}$.  

{Equation} \cref{eq:vmsrom} is referred to as VMS-ROM and can be written as:
\begin{align}
    B \frac{d \uu_{\tr}}{dt} & = -\mC (\bar{\uu}_{\tr}) \bar{\uu}_{\tr} -
    \frac{1}{\rm Re} A \bar{\uu}_{\tr} + \utau_{\text{VMS}}, 
    \label{eq:vmsrom_ode}
\end{align}
where 
\begin{equation}
    \tau_{\text{VMS},i} \equiv
    \biggl( 
        (\bu_R \cdot \nabla) \bu_R, \bphi_{i}
    \biggr) - \biggl( 
        (\bu_r \cdot \nabla) \bu_r, \bphi_{i}
    \biggr), \qquad \forall~i=1,\ldots,r. \label{eq:truth-vmsrom-closure}
\end{equation}

\begin{remark}
    If we drop the VMS-ROM closure term, we are left with the G-ROM (\ref{eq:gromu}).
\end{remark}

\begin{remark}
    The VMS-ROM closure term is a \textit{correction term} in the higher-dimensional space $\bX_1 \bigoplus \bX_2$.
\end{remark}

The VMS-ROM closure term $\utau_{\text{VMS}}$ is essential for the accuracy of \cref{eq:vmsrom_ode}. However, the unresolved component of $\bu_\tR$, $\bu'_\tr$, is not available during the online stage. 
To overcome this issue, for each $i=1,\ldots,r$, $\tau_{\text{VMS},i}$ is approximated using a \textit{generic} function $g_i(\uu_\tr)$ 
as follows:
\begin{align}
    \tau_{\text{VMS},i} \equiv \biggl( 
        (\bu_\tR \cdot \nabla) \bu_\tR, \bphi_{i}
    \biggr)  
    - \biggl( 
        (\bu_\tr \cdot \nabla) \bu_\tr, \bphi_{i}
    \biggr)  \approx 
    g_i(\uu_\tr).
\end{align}
In the {VMS-ROM} training, for each $i=1,\ldots,r$, the function $g_i$ is determined by minimizing the 
{mean squared error (MSE)}
between the VMS-ROM closure $\tau_i$ and the reconstruction made by $g_i$, 
\begin{equation}
    \text{MSE}_{\text{tr}, i} := \frac{1}{N_{\rm tr}} \sum_{j=1}^{N_{\rm tr}} (\tau^j_i - 
    g_i(\uhu^j_\tr))^2. 
    \label{eq: mse_tr}
\end{equation}
We note that $N_{\rm tr}$ is the number of training time samples, {which,} 
{in our numerical experiments,} is selected to be the same as the number of snapshots that are used to build the POD basis functions. 
The VMS-ROM closure term $\tau^j_i$ is defined to be
\begin{equation}
  \tau^j_i \equiv
      \biggl( 
       (\mP_\tR\bu^j \cdot \nabla) \mP_\tR \bu^j, \bphi_{i}
   \biggr)  
   - \biggl( 
       (\mP_\tr \bu^j \cdot \nabla) \mP_\tr\bu^j, \bphi_{i}
   \biggr), 
  \label{eq:vmsrom-closure}
\end{equation}
where $\mP_\tr$ is 
{the} operator that projects the FOM solution $\bu \in \bX^\cN$
onto the $r$-dimensional reduced space $
{\bX_r}$, and $\bu^j$ is the FOM
solution at time $t^j$. The 
vector {of ROM coefficients} $\uhu_\tr$ is obtained by projecting the FOM solution onto the $r$-dimensional reduced space, that is, $\widehat{u}^j_{r,i} = (\bu^j,\bphi_i)$ for all $i=1,\ldots,r$. We note that for
$i=1,\ldots,r$, the function $g_i$ is found independently.
\begin{remark}
   By {a slight} abuse of notation, we refer to both (\ref{eq:truth-vmsrom-closure}) and (\ref{eq:vmsrom-closure}) as the VMS-ROM closure term. We note that (\ref{eq:vmsrom-closure}) is used to train the function $g_i$ instead of (\ref{eq:truth-vmsrom-closure}) because 
   {the latter} requires running {the $R$-dimensional} G-ROM, which is not practical. Moreover,
   the stiffness terms are ignored in both (\ref{eq:truth-vmsrom-closure}) and (\ref{eq:vmsrom-closure}) because 
   {these terms are} usually small {in the convection-dominated regime} due to the $1/\rm Re$ scaling {(see the careful discussion of the commutation error in \cite{koc2019commutation})}. 
\end{remark}

\begin{figure}[!ht]
    \centering     
    \includegraphics[width=0.8\columnwidth]{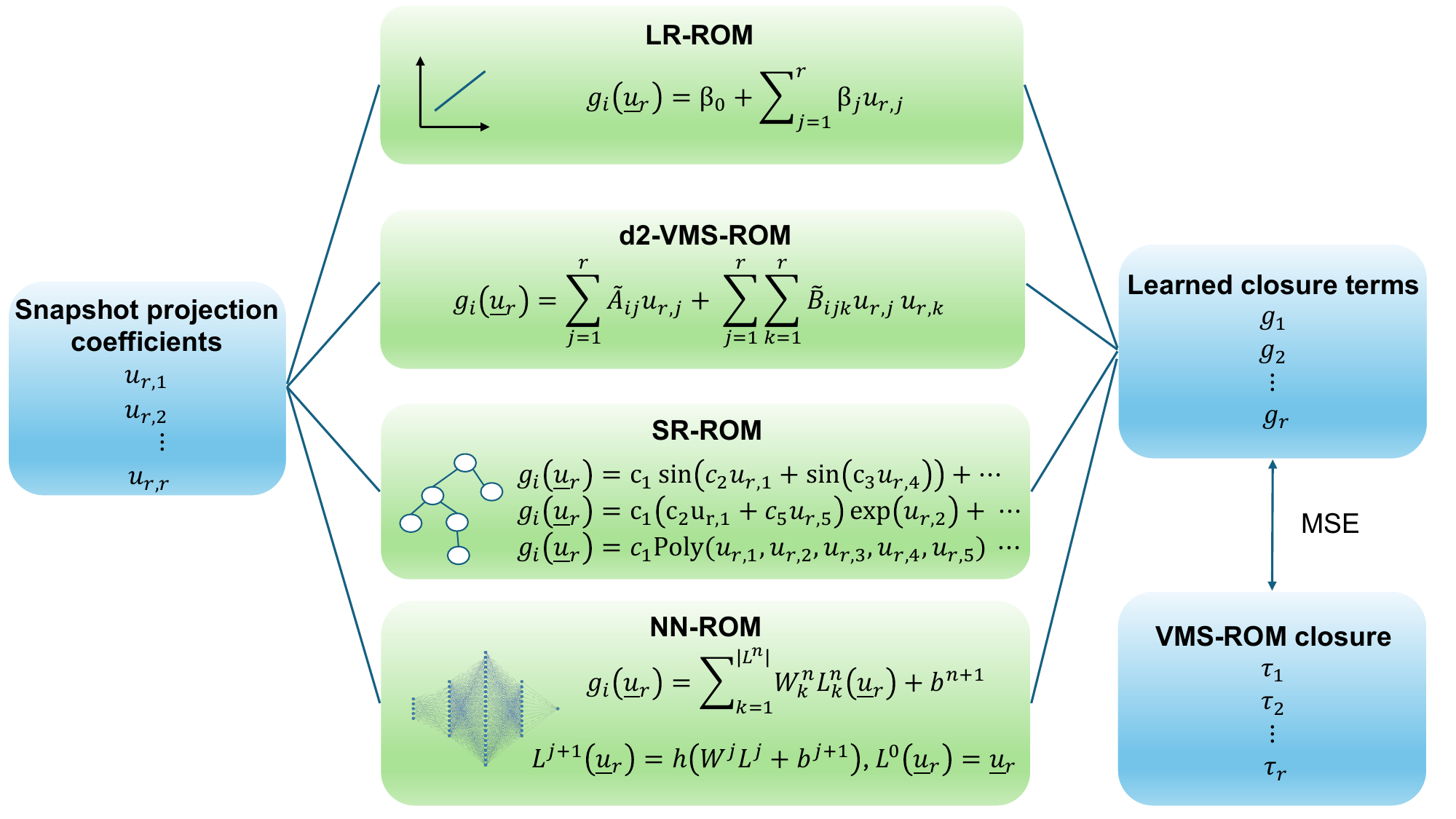}
    \caption{Schematic of the 
    {ML-VMS-ROMs} considered
    in this paper.
    }
    \label{fig:ml_vms_rom_schematic}
\end{figure}

In the 
{remainder of this section}, we describe {the} four 
ML techniques that we {employ} to construct the function $g_i$ in the ML-VMS-ROMs:
the linear regression VMS-ROM (outlined in Section~\ref{subsec: lr-rom});
the data-driven VMS-ROM (outlined in Section~\ref{subsec:d2-vms-rom});
the novel symbolic regression VMS-ROM (outlined in Section~\ref{subsec:sr-rom}); and
the neural network VMS-ROM (outlined in Section~\ref{subsec:nn-rom}). A schematic of the four different ROM closure modeling strategies that are 
used to construct the ML-VMS-ROMs
is depicted in Figure~\ref{fig:ml_vms_rom_schematic}.

\subsubsection{VMS-ROM with Linear Regression Closure Models (LR-ROM)} 
\label{subsec: lr-rom}
Linear regression (LR) predicts a target (dependent) variable $y$ (in our case, the $i$-th component of the VMS ROM closure, $\tau_i$) based on one or more input (independent) variables $x_j$ (in our case, the $j$-th ROM coefficient, $u_{r,j}$) by assuming a linear dependency between them. Mathematically, for each $i=1,\ldots, r$, the function $g_i$ is postulated as
\begin{eqnarray}
    g_i(\uu_\tr) = \beta_0 + \sum_{j=1}^r \beta_j u_{r,j},
\end{eqnarray}
where $r$ is the dimension of the reduced space and $(\beta_{0}, \dots, \beta_r)$ are the learnable weights. 
The optimal weights are found by solving the linear system given by the 
normal equations 
(see equation (2.5) in \cite{hastie2009elements}). To penalize the complexity of the model and avoid overfitting, 
{the} regularization term 
$\alpha \sum_{j=1}^r \beta_{j}^2$ is added to \eqref{eq: mse_tr} to favor the solution with weights of smaller $L^2$ norm. 

In this work, the \texttt{Ridge} class of the \texttt{scikit-learn} library is used for training the LR-ROM. The default settings are used, apart for the hyperparameter $\alpha$, which is tuned through a hold-out validation based on the error metric (\ref{eq: mse_val}). A more detailed discussion on LR can be found in \cite{gareth2013introduction}. The advantages of LR stem from its interpretability and simplicity. 
However, linear models may not always effectively describe data, and they are particularly sensitive to outliers.

\subsubsection{Data-Driven Variational Multiscale ROMs (d2-VMS-ROM)}
\label{subsec:d2-vms-rom}

In \cite{mou2021data}, the authors proposed a data-driven approach to
approximate the VMS-ROM closure. In particular, the structure of function $g_i$ is assumed to be similar to the structure of the 
{G-ROM}
function and is
postulated as 
\begin{eqnarray}
    g_i(\uu_r) = \sum_{j=1}^r \tA_{ij} u_{r,j} + \sum_{j=1}^r \sum_{k=1}^r \tB_{ijk} u_{r,j}u_{r,k}.
\end{eqnarray}
The entries in the $r \times r$ matrix $\tA$ and the $r \times r \times r$ tensor $\tB$ are determined by minimizing the training MSE (\ref{eq: mse_tr}). We note that the minimization problem is equivalent to the following least squares problem
\begin{align}
  \min_{\ux} \|G\ux -\uf\|^2, \label{eq:lstq}
\end{align}
where $G \in \mathbb{R}^{N_{\rm tr} r \times ({r^2+r^3})}$ is a matrix whose entries are determined by $\uhu^j_\tr$, for $j=1,\ldots,N_{\rm tr}$, $\uf \in \mathbb{R}^{N_{\rm tr} r \times 1}$ is a vector whose entries are determined by $\utau^j$ for $j=1,\ldots,N_{\rm tr}$, and $\ux \in \mathbb{R}^{(r^2+r^3)\times 1}$ is a vector whose entries are determined by $\tA$ and $\tB$.
To alleviate the ill-conditioning of the least squares problem (\ref{eq:lstq}) \cite{peherstorfer2016data}, we use the truncated SVD as a regularization method \cite{mou2021data}.

\subsubsection{VMS-ROM with Symbolic 
Regression Closure Models (SR-ROM)}
\label{subsec:sr-rom}

Symbolic regression (SR) is a type of regression analysis that seeks to find mathematical expressions that accurately describe the relationship between the features and the labels of a dataset \cite{kronberger2024symbolic}. Unlike traditional regression methods, which rely on predefined model structures, SR
identifies both the structure of the model and its parameters. 

Here, we use SR to search for the function $g_i$ in (3.14) that solves the following optimization problem
\begin{equation}
    g_i = \argmax_{f \in \mathcal{F}} \left[ 1 - \frac{\sum_{j=1}^{N_{\rm tr}} (\tau^j_i - 
    f(\uhu^j_\tr))^2}{\sum_{j=1}^{N_{\rm tr}} (\tau_i^j - \overline{\tau_i})^2} \right], \label{eq: sr_train}
\end{equation}
where $\overline{\tau_i}$ is the mean of ${\tau}_i$ over $\Ntr$ samples and $\mathcal{F}$ is a set of functions generated by a \textsl{primitive set}, that is, a collection of mathematical operations, such as $\{+,-,*,\sin, \exp\}$, and a \textsl{terminal set}, made of the independent variables and some constants. We note that maximizing the objective function in \eqref{eq: sr_train} (called \textsl{fitness}) is equivalent to minimizing the MSE in \eqref{eq: mse_tr}. Apart from dealing with algebraic expressions, SR has recently been extended to field problems by including discrete differential operators \cite{manti2024}.

\begin{figure}[!ht]
    \centering
     \includegraphics[width=0.85\textwidth]{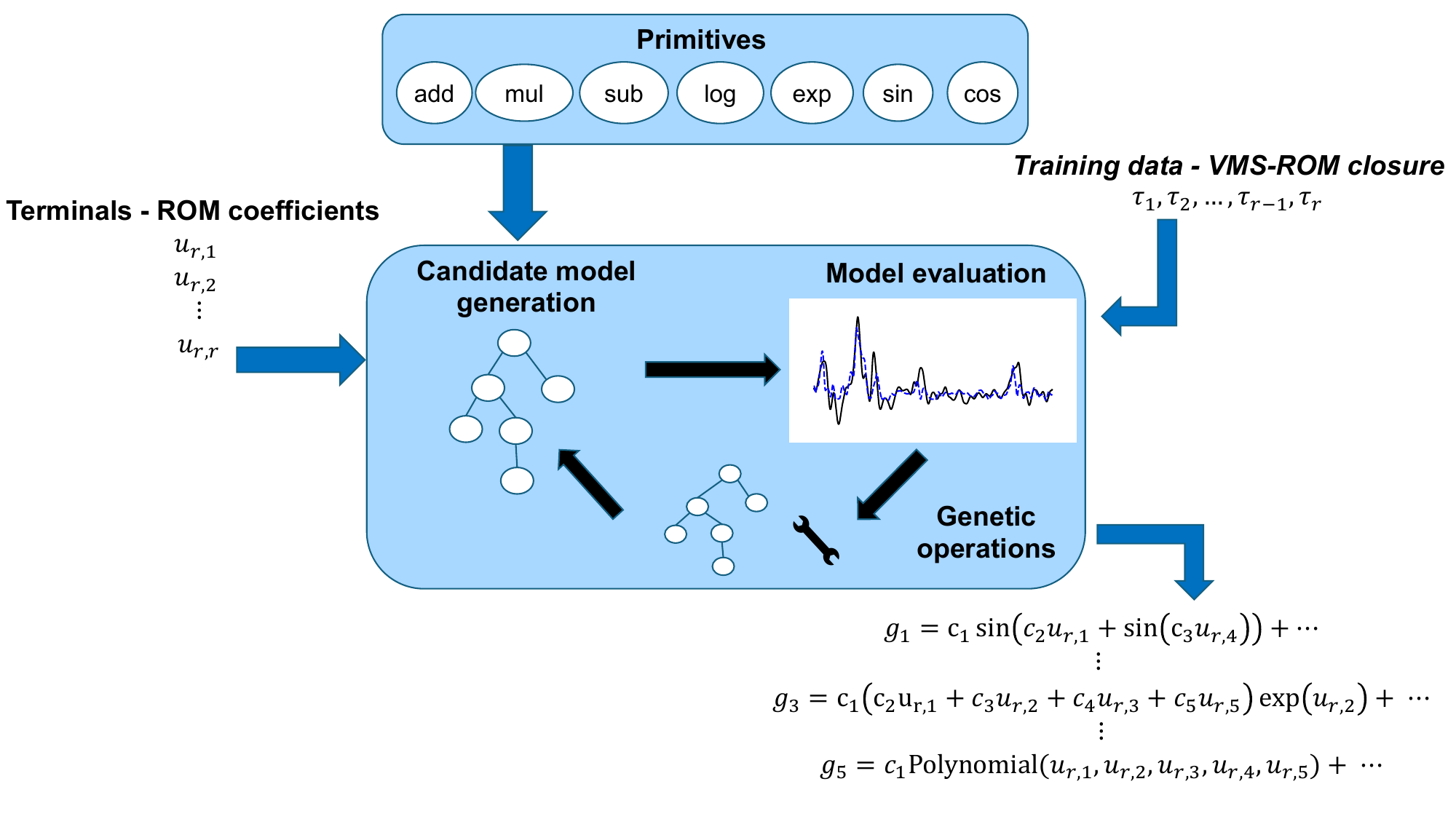}
    \caption{Schematic of the symbolic regression approach {used to model ROM closures (SR-ROM).}
    }
    \label{fig:sr_ls}
\end{figure}

A common approach to SR (also adopted in this work) is genetic programming (GP) \cite{koza1994genetic, o2009riccardo}. GP is a gradient-free strategy that explores a space of mathematical expressions by iteratively evolving a population of candidate models through genetic operations. A candidate model, is also known as \textsl{individual}, can be represented by a tree, where each node is either a \textsl{primitive} or a \textsl{terminal} (\textsl{i.e.}, a variable or a constant), 
and is characterized by its \textsl{fitness}, as defined above.  GP-based SR consists of the following steps:
\begin{enumerate}
    \item An initial population of individuals is randomly generated, based on the chosen primitive set;
    \item For each individual, the constants are optimized via the Levenberg–Marquardt method to maximize the fitness function;
    \item A pool of individuals are \textsl{selected} according to fitness-based tournaments. Specifically, for each tournament a given number (tournament size) of individuals are sampled from the population and the ones with the highest fitness are selected for the next step;
    \item The selected expressions are modified using \textsl{crossover} (combining parts of two parent expressions to create new offspring) and \textsl{mutation} (randomly altering parts of an expression);
    \item Steps 2–4 are repeated for multiple \textsl{generations}, until a stopping condition is satisfied or a maximum number of generations is reached.
\end{enumerate}
A more detailed review of SR methods can be found in  \cite{koza1994genetic, o2009riccardo}. A schematic of the described algorithm is reported in Figure~\ref{fig:sr_ls}.

In this paper, we used the SR implementation of the \texttt{pyoperon} \cite{Burlacu2020} library, because it achieves a good compromise between test performance, final model size, and training time (see, \textsl{e.g.}, Figure~1 in \cite{la2021contemporary}). We fixed the following hyperparameters for every run:
\begin{itemize}
    \item tournament size equal to $5$;
    \item number of iterations of the constant optimizer equal to 10;
    \item population size equal to 1000;
    \item maximum depth of the trees equal to 5.
\end{itemize}
The maximum length of the individuals (\textsl{i.e.}, the number of nodes), the number of generations, and the primitive set were tuned through a grid search based on a hold-out validation. A description of the datasets 
used for training, validation, and testing can be found in Section~\ref{sec:nr}. 

\subsubsection{VMS-ROM with Neural Network Closure Models (NN-ROM)}
\label{subsec:nn-rom}

Among the different architectures of neural networks (NNs), in this work, we used a multi-layer perceptron (MLP) to represent the mapping between the vector of ROM coefficients $\uu_r$ and the approximated closure term. The MLP can be represented mathematically as follows:
\begin{align}
    &g_i(\uu_r) := \sum_{k=1}^{|L^n|} W^n_{k} L^{n}_k(\uu_r)  + b^{n+1},\\
    & {L}^{j+1}(\uu_r) := h({W}^j {L}^j + {b}^{j+1}), \quad 0 \leq j \leq n-1, \\
    & {L}^0(\uu_r) := \uu_r,
\end{align}
where $n$ is the number of the (hidden) layers, ${L}^j$ is the output of the $j$-th layer, $|L^n|$ is the number of nodes of the $n$-th layer, $h$ is the activation function, ${W}^j$ is the weight matrix for the connections between the $j$-th and the $(j+1)$-th layers, and ${b}^j$ is the bias term for the outputs of the $j$-th layer. 

The training of the NN amounts to adjusting its weights to minimize a loss function using \textsl{backpropagation} \cite{rumelhart1986learning}. Here, as a loss function, we chose the sum of the MSE \eqref{eq: mse_tr} and the $L^2$ regularization term, which penalizes the sum of the squares of the weights. A more detailed introduction to NNs can be found in \cite{Goodfellow-et-al-2016, gareth2013introduction}. 

In this work, we used the \texttt{skorch} \cite{tietz2017skorch} library for training the MLP with the following settings:
\begin{itemize}
    \item \texttt{LeakyReLU} as the activation function;
    \item Adam optimizer;
    \item 100 training epochs;
    \item batch size equal to 512.
\end{itemize}
We tuned the values of the learning rate, the number of units per layer, the $L^2$ regularization parameter, and the dropout rate through a grid search based on a hold-out validation. Details on the datasets used for training, validation, and testing can be found in Section~\ref{sec:nr}. 

\section{Numerical Results}
\label{sec:nr}
In this section, we assess the four ML-VMS-ROMs introduced in Section~\ref{subsection:d2-vms-roms} in the numerical simulation of two test problems: (i) the 2D flow past a cylinder, and (ii) the 2D lid-driven cavity. We emphasize that, for both test problems, we consider the under-resolved regime, i.e., the case when the number of ROM basis functions is not large enough to capture the underlying dynamics.
We also note that, for the 2D lid-driven cavity, we consider the convection-dominated regime. 
Both the under-resolved regime and the convection-dominated regime are relevant to the realistic engineering and geophysical settings (\textit{e.g.}, the numerical simulation of turbulent flows).  

We remark that each model depends on hyperparameters, and we list the considered hyperparameters for each model below. 
In each model, the function $g_i$ is determined by minimizing the training MSE, $\text{MSE}_{\text{tr}, i}$ (\ref{eq: mse_tr}). 
\begin{itemize}
    \item For the LR-ROM, we select the $L^2$ regularization parameter 
    from the set
    $\{10^i, 2 \times 10^i, \dots, 9 \times 10^i\}_{i=0,\dots,5}$.
    \item For the d2-VMS-ROM, we select the SVD rank 
    from the set
    $\{1,2,\ldots,\text{Rank($G$)}\}$, where $G$ is the matrix in the least squares problem (\ref{eq:lstq}), whose entries are determined by the ROM coefficients of the projected snapshots, $\uhu_\tr$.
    \item For the SR-ROM, we select the maximum length 
    from the set
    $\{5,10,\dots,45,50\}$, the number of generations 
    from the set
    $\{10, 25, 50, 75, 100\}$, and the primitive set 
    from the set \{\{\texttt{add,sub,mul,constant,variable}\}, \{\texttt{add,sub,mul,sin,constant,variable}\}, \{\texttt{add,sub,mul,exp,sin,constant,variable}\}, \{\texttt{add,sub,}  \texttt{mul,sin,cos,constant,variable}\}, \{\texttt{add,sub,mul,exp,sin,cos,constant,variable}\}, \{\texttt{add,sub,} \texttt{mul,exp,sin,cos,square,log,} \texttt{constant,variable}\} \}\};
    \item For the NN-ROM, we select the learning rate 
    from the set
    $\{10^{-4}, 10^{-3}, 10^{-2}\}$, the number of nodes per layer 
    from the set
    $\{[64, 128, 256, 512,256, 128, 64], [64, 128, 256, 128, 64], [64, 128, 64]\}$, the $L^2$ regularization parameter 
    from the set
    $\{10^{-5}, 10^{-4}, 10^{-3}\}$, and the dropout rate 
    from the set
    $\{0.3,0.4,0.5\}$. 
\end{itemize}
For each model and each hyperparameter, the selected function $g_i$ is then integrated with the VMS-ROM, and its performance is further evaluated in the validation interval. The optimal hyperparameters are selected to be 
those that maximize the coefficient of determination
\begin{equation}
    R^2_{\rm val} := 1 - \frac{\text{MSE}_{\text{val}}}{\frac{1}{N_{\rm val}} \sum_{i=1}^{N_{\rm val}} (E^i_{\rm FOM} - \overline{E_{\rm FOM}})^2},
\end{equation}
where the validation error, $\text{MSE}_{\text{val}}$, is defined to be the MSE between the FOM kinetic energy and the ROM 
kinetic energy,
\begin{equation}
    \text{MSE}_{\text{val}} = \frac{1}{N_{\rm val}} \sum_{i=1}^{N_{\rm val}} (E^i_{\rm FOM} - E^i_{\rm ROM})^2. \label{eq: mse_val}
\end{equation}
The kinetic energy $E$ at time $t^j$ is defined to be
\begin{equation}
    E(t^j) = E^j = \frac{1}{2} \left(\int_{\Omega} \bu^j \cdot \bu^j~d\Omega \right)^{\frac{1}{2}}, 
\end{equation}
and the averaged kinetic energy over the $N_{\rm val}$ samples in the validation set is defined to be
\begin{equation}
    \overline{E_{\rm FOM}} = \frac{1}{N_{\rm val}} \sum_{i=1}^{N_{\rm val}} E^i_{\rm FOM}.
\end{equation}

Once the optimal parameters are selected,  the resulting ML-VMS-ROM is assessed in a testing interval with the metric
\begin{equation}
    \text{rMSE}_{\text{test}} := \frac{1}{N_{\rm test}} \sum_{i=1}^{N_{\rm test}} \frac{(E^i_{\rm FOM} - E^i_{\rm ROM})^2}{(E^i_{\rm FOM})^2}. \label{eq: mse_test}
\end{equation}

\subsection{2D Flow Past a Cylinder}
\label{subsec: 2dfpc}

\begin{figure}[!ht]
    \centering
     \includegraphics[width=0.95\textwidth]{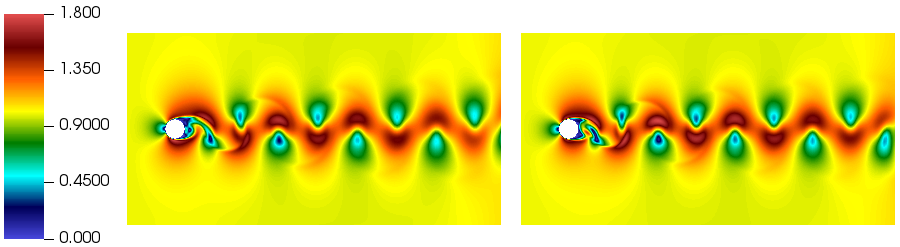}
    \caption{2D flow past a cylinder. Velocity magnitude snapshots for $\rm Re=400$
    ({left}) and $\rm Re=500$
    ({right}).} 
    \label{fig:2dfpc_fom_snap}
\end{figure}

Our first test problem is the 2D flow past a cylinder, which is a canonical test
case for ROMs. The computational domain is $\Omega = [-2.5 D:17 D] \times [-5 D:5 D]$,
where $D$ is the cylinder diameter, and the cylinder is centered at $[0,0]$. 
We consider two Reynolds numbers, namely, $\rm Re=400$ and $\rm Re=500$. For each Reynolds number, we collect $\Ntr=2001$ snapshots, $\{\bu^k := \bu(\bx,t^k)-\bphi_0\}^{\Ntr}_{k=1}$, over the time interval $[500,~520]$ (measured in
convective time units, $D/U$, where $U$ is the free-stream
velocity). These snapshots are taken after the von Karman vortex street has developed, with sampling time $\Delta t_s=0.01$. The zeroth mode, $\bphi_0$, is set to be the FOM velocity
field at $t=500$. 
Snapshots of the velocity magnitude for $\rm Re=400$ and $\rm Re = 500$ are displayed in Fig.~\ref{fig:2dfpc_fom_snap}.

The snapshot data are then used to construct the reduced basis functions $\{\bphi_i\}^r_{i=1}$ by using POD. The ROM coefficients,
$\{\uhu^j_{\tr}\}^{\Ntr}_{j=1}$, are computed by projecting the snapshots onto the $R$-dimensional reduced space, where $R$ is the rank of the snapshot matrix. The VMS-ROM closure terms $\{\utau^j\}^{\Ntr}_{j=1}$ are then computed  using~\eqref{eq:vmsrom-closure}. 
For each ML-VMS-ROM, the function $g_i$ is found by minimizing the training MSE \eqref{eq: mse_tr}. 
We select the optimal values for the hyperparameters for model training such that they minimize the validation error~\eqref{eq: mse_val}. To evaluate this error, we further run the FOM simulation and monitor the FOM kinetic energy over the time interval $[520,~540]$. The trained function $g_i$ is then used within the VMS-ROM to perform a simulation over the time interval $[500,~540]$. Finally, we assess the generalization error in \eqref{eq: mse_test} of all the ROMs in an extrapolation region, that is, in the time interval $[540,~600]$.

In Table~\ref{tab:2dcyl_test_error}, the relative MSE ($\text{rMSE}_{\text{test}}$) 
between the FOM kinetic energy and the predicted kinetic energy during the time interval $[540, 600]$ (extrapolation region) for all the ML-VMS-ROMs is reported with respect to the reduced space dimension $r=2,3,4,5$. 
Because the models trained using neural networks and symbolic regression rely on stochastic optimization algorithms, the $\text{rMSE}_{\text{test}}$ for both the NN-ROM and the SR-ROM are computed by averaging the results over 5 runs corresponding to different seeds.
\begin{table}[!ht]
    \centering
    \caption{2D flow past a cylinder. Relative MSE ($\text{rMSE}_{\text{test}}$) between the FOM kinetic energy and the kinetic energy predicted in the extrapolation region (\textsl{i.e.}, in the time interval $[540,600]$) by the 
    ML-VMS-ROMs,
    for 
    different values of the reduced space dimension $r=2,3,4,5$. For the SR-ROM and the NN-ROM, the mean and the standard deviation over $5$ independent runs (different seeds) are reported. 
    In bold, the best performing model for each $r$ value {is listed}.}
    \begin{tabular}{cccccc}
        \toprule
         Re &  $r$ & LR-ROM & d2-VMS-ROM & SR-ROM & NN-ROM\\
         \midrule
         \multirow{4}{*}{400} & 2 &$1.4 \times 10^{-4}$ & $5.3 \times 10^{-6}$ & \boldmath $(7.1 \pm 1.1) \times 10^{-7}$ & $(5.3 \pm 3.7) \times 10^{-6}$\\
         & 3 & $3.3 \times 10^{-4}$ & nan & \boldmath $(1.5 \pm 2.0) \times 10^{-5}$ & $(3.1 \pm 3.1) \times 10^{-4}$\\
         & 4 & $5.5 \times 10^{-5}$ & \boldmath $7.5 \times 10^{-7}$ & \boldmath $(6.7 \pm 1.0) \times 10^{-7}$ & $(4.4 \pm 2.7) \times 10^{-6}$\\
         & 5 & $4.6 \times 10^{-6}$ & $2.0 \times 10^{-5}$ & \boldmath $(9.5 \pm 8.9) \times 10^{-7}$ & \boldmath $(10.0 \pm 6.7) \times 10^{-7}$\\
         \midrule
         \multirow{4}{*}{500} & 2 & $8.1 \times 10^{-6}$ & $3.2 \times 10^{-5}$ & $(1.7 \pm 1.9) \times 10^{-5}$ & \boldmath $(4.5 \pm 1.6) \times 10^{-6}$\\
         & 3 & $6.6 \times 10^{-5}$ & $3.7 \times 10^{-5}$ & \boldmath $(1.8 \pm 1.0) \times 10^{-6}$ & \boldmath $(2.7 \pm 1.6) \times 10^{-6}$\\
         & 4 & $2.2 \times 10^{-5}$ & \boldmath $3.2 \times 10^{-6}$ & \boldmath$(7.2 \pm 5.1) \times 10^{-6}$ & \boldmath$(10.4 \pm 7.8) \times 10^{-6}$\\
         & 5 & $3.1 \times 10^{-4}$ & $1.9 \times 10^{-4}$ & \boldmath $(1.0 \pm 1.3) \times 10^{-6}$ & $(4.7 \pm 2.3) \times 10^{-5}$ \\
         \bottomrule
    \end{tabular}
    \label{tab:2dcyl_test_error}
\end{table}
The results show that the SR-ROM is the most robust model with respect to $r$ for both $\rm Re$, followed by the NN-ROM, the d2-VMS-ROM, and the LR-ROM. 
In particular, we note that for $\rm{Re} = 400$ and $r=2$ and $r=3$ the SR-ROM $\text{rMSE}_{\text{test}}$ values
in the extrapolation region are an order of magnitude lower than the $\text{rMSE}_{\text{test}}$  values
of the other ROMs. The same happens for $\rm{Re} = 500$ and $r=5$.

We further measure the closeness between the VMS-ROM closure terms $\{\utau^j\}^{\Ntr}_{j=1}$ and the reconstructed closure terms for each model. In particular, the MSE for each $\utau$ component, $\text{MSE}_{\text{tr}, i}$~\eqref{eq: mse_tr} is computed, and the mean MSE, which is defined as
\begin{equation}
    \overline{\text{MSE}_{\text{tr}}} := \frac{1}{r} \sum_{i=1}^r \text{MSE}_{\text{tr}, i}, \label{eq: mean_mse_tr}
\end{equation}
is reported in Table~\ref{tab:2dcyl_tau_reconstruction_error}.
We observe that SR-ROM is the most robust method also in this case. We also note that the NN-ROM has a much larger MSE than the SR-ROM even though the corresponding test errors on predicted energy are comparable (Table~\ref{tab:2dcyl_test_error}). 
For example, for $\rm Re = 500$, 
the training MSE of NN-ROM with $r=3$ is roughly two
orders of magnitude larger than the training MSE of SR-ROM. However, the performance of the two models in the extrapolation region is comparable.
There are two reasons for this behavior:
\begin{enumerate}
    \item {We are using different metrics for training (closure error) and validation/testing (kinetic energy error)}.
    \item In the NN-ROM and LR-ROM,  we consider $L^2$ regularization to prevent models from overfitting. This results 
    in a model with smaller weights and could lead to a larger training MSE but better generalizability (smaller $\text{rMSE}_{\text{test}}$ in the {extrapolation} region).
\end{enumerate}
{Overall, our results demonstrate that, SR-ROM outperforms the other ML-VMS-ROM models in terms of generalizability.}

\begin{table}[!tb]
    \centering
    \caption{2D flow past a cylinder. Mean squared error ({$\overline{\text{MSE}_{\text{tr}}}$}) between the VMS-ROM closure and the closures given by the 
    ML-VMS-ROMs 
    for $\rm Re=400,500$. 
    In bold, the best performing model for each case is listed.}
    \begin{tabular}{cccccc}
        \toprule
         Re &  $r$ & LR-ROM & d2-VMS-ROM & SR-ROM & NN-ROM\\
         \midrule
         \multirow{4}{*}{400} & 2 & $2.3 \times 10^{-3}$ & $5.2 \times 10^{-4}$  & \boldmath $(6.9 \pm 12.5) \times 10^{-7}$ & $(8.5 \pm 4.6) \times 10^{-4}$\\ 
         & 3 & $5.9 \times 10^{-3}$ & $5 \times 10^{-4}$ & \boldmath $(1.3 \pm 2.6) \times 10^{-4}$ & $(1.7 \pm 1.1) \times 10^{-3}$\\
         & 4 & $3.1 \times 10^{-3}$ & $6.1 \times 10^{-4}$ &  \boldmath $(1.2 \pm 1.5) \times 10^{-4}$ & $(1.2 \pm 0.5) \times 10^{-3}$\\
         & 5 & $1.4 \times 10^{-3}$ & \boldmath $5.6 \times 10^{-6}$ & $(5.7 \pm 4.9) \times 10^{-5}$ & $(1.9 \pm 0.4) \times 10^{-4}$\\
         \midrule
         \multirow{4}{*}{500} & 2 & $4.7 \times 10^{-4}$ & $1.2 \times 10^{-3}$ & \boldmath $(5.5 \pm 3.3) \times 10^{-5}$ & $(8.7 \pm 4.3) \times 10^{-5}$ \\ 
         & 3 & $5.9 \times 10^{-3}$ & $3.5 \times 10^{-4}$ & \boldmath $(7.5 \pm 12) \times 10^{-6}$ & $(5.5 \pm 3.6) \times 10^{-4}$\\
         & 4 & $3.5 \times 10^{-3}$ & $4.1 \times 10^{-3}$ & \boldmath $(1.4 \pm 2.7) \times 10^{-4}$ & $(4.9 \pm 2.7) \times 10^{-4}$\\
         & 5 & $1.0 \times 10^{-3}$ & $4 \times 10^{-4}$ & \boldmath $(5.2 \pm 3.0) \times 10^{-5}$ & $(8 \pm 6.1) \times 10^{-4}$\\
         \bottomrule
    \end{tabular}
    
    \label{tab:2dcyl_tau_reconstruction_error}
\end{table}

\begin{figure}[!ht]
    \centering
    \includegraphics[scale=0.52]{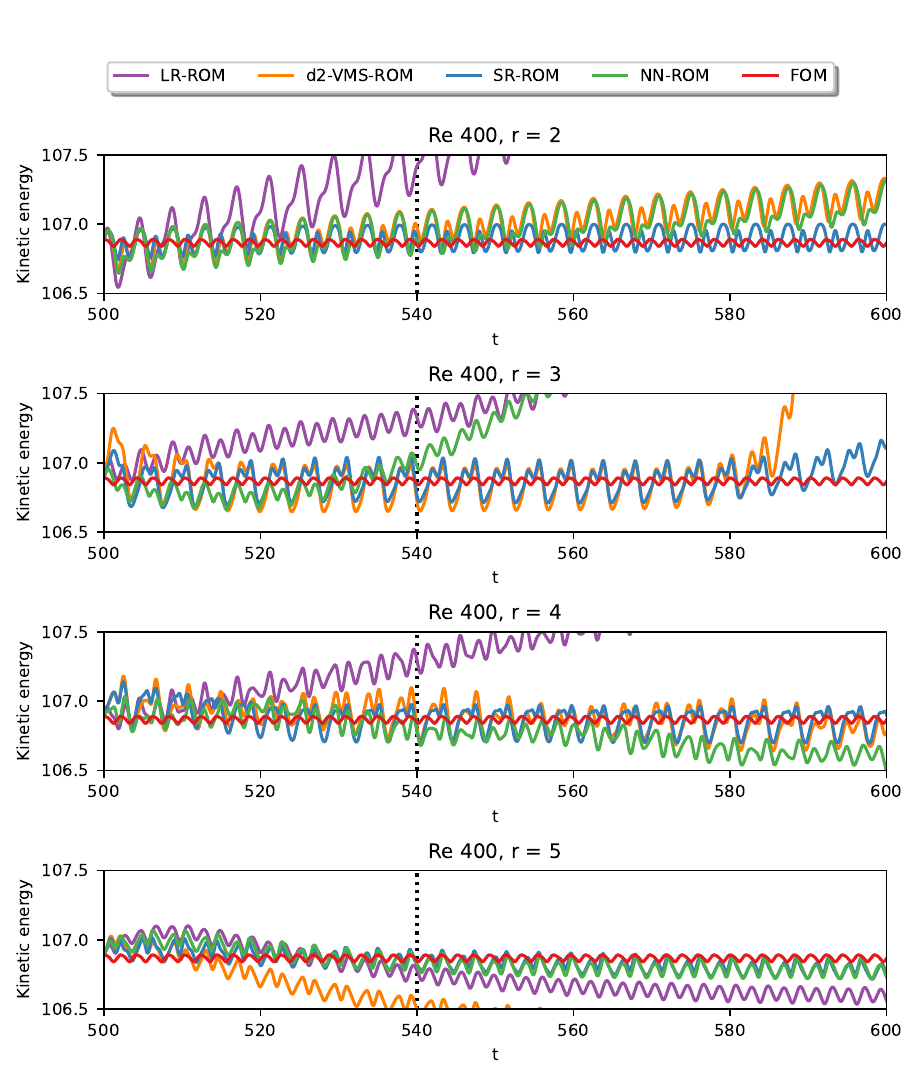}
     \includegraphics[scale=0.52]{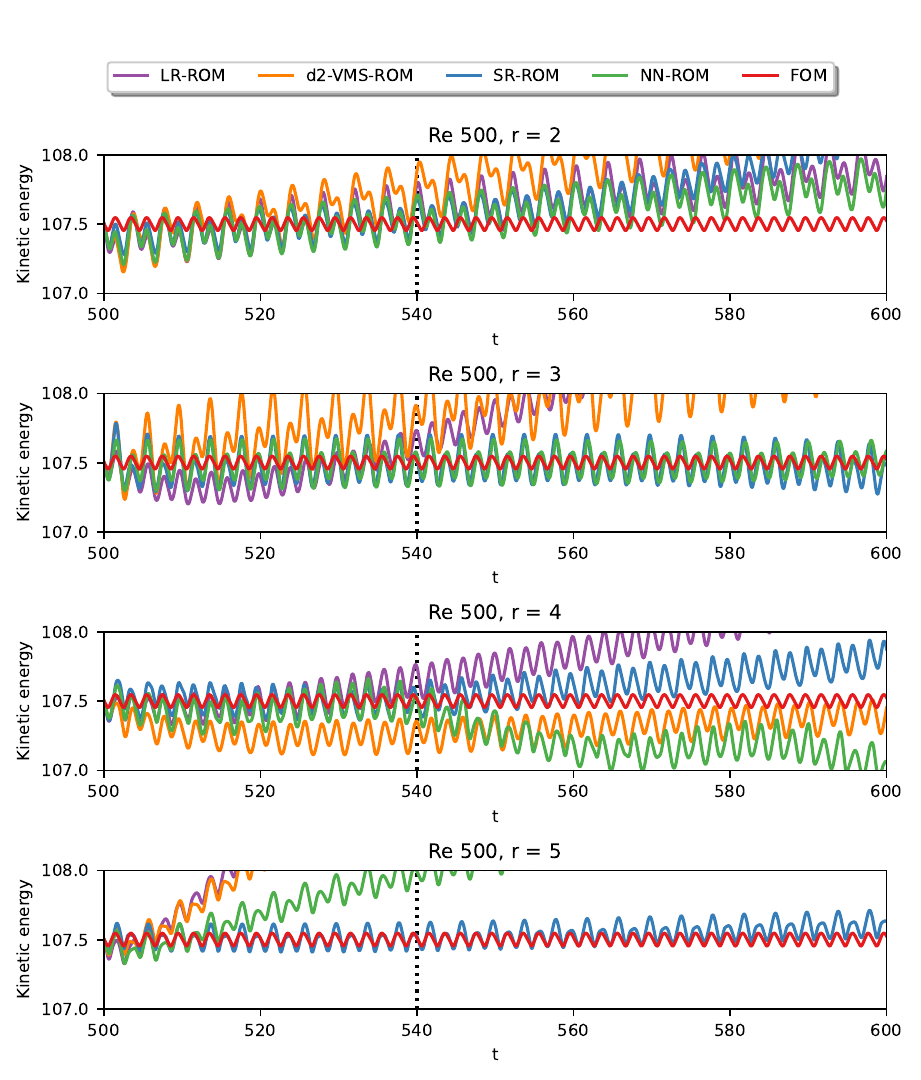}
    \caption{2D flow past a cylinder at $\rm Re=400$ and $\rm Re=500$. Comparison of the kinetic energy of the FOM and the ROMs for the reduced space dimension $r=2,3,4,5$.}
    \label{fig:energies_2dcyl}
\end{figure}

In Fig.~\ref{fig:energies_2dcyl},
we show the kinetic energy behavior of the ROMs and the FOM in the time interval $[500, 600]$ with reduced space dimension $r=2,3,4,5$ for $\rm Re=400$ and $\rm Re=500$. 
For $\rm Re=400$, we observe that the LR-ROM energy 
deviates away from the FOM energy for all $r$ values considered. 
For the d2-VMS-ROM and the NN-ROM, although the energy is better than the LR-ROM energy, it could still deviate away from the FOM energy; see, for example, the $r=3$ case. The SR-ROM is the most stable compared to the other three ML-VMS-ROMs, and it is able to predict energy that is comparable to the FOM energy.
For $\rm Re=500$, we again observe that the 
LR-ROM energy deviates away from the FOM energy for all $r$ values considered. 
For the d2-VMS-ROM and the NN-ROM, although the energy is better than the 
LR-ROM energy, it could still deviate away from the FOM energy; see, for example, the $r=5$ case. The SR-ROM is the most stable with respect to different $r$ values.

 \begin{table}[!ht]
    \caption{2D flow past a cylinder. 
    Expressions of the components of the SR-ROM closure model that attains the minimum $\text{rMSE}_{\text{test}}$ in the extrapolation region reported in Table~\ref{tab:2dcyl_test_error}.}  
	\begin{center}
		\begin{tabular}{c c c}
			\toprule
			Re & Component &  Expression\\
			\midrule
            \multirow{4}{*}{400} &  $g_1$ &\cellcolor{lightgray} \small $(-0.039u_1 + 0.002u_2 - 0.034)\sin(1.719u_2 - 3.789) + 0.000$\\
            &$g_2$ & \small $(-0.077u_1 + 0.031u_4 - 0.073)\sin(1.062u_1 - 1.804u_2) + 0.020$\\
            &$g_3$ & \cellcolor{lightgray} \small $-0.384u_1 + 0.384u_3 + 0.039u_4 + 0.971\sin(1.868u_1 + 0.123) - 0.012$\\
            &$g_4$ & \small $(7.533u_2 - 1.921u_4)(0.706u_2 + 0.046u_4 - \sin(0.846u_2)) - 0.008$\\
            \hline
            \multirow{5}{*}{500} & $g_1$ &\cellcolor{lightgray} \small $0.015u_1 + 0.010u_4 + \sin(0.191u_1 - 0.100u_5 + \sin(0.064u_4) - 1.569) + 0.979$\\
            &$g_2$ & \small $1.004(0.264u_2u_3 - \exp(0.282u_2))(0.007u_1 - 0.005u_2 - 0.014u_3 + 0.006u_4) - 0.002$\\
            &$g_3$ & \cellcolor{lightgray} \small $-0.999(0.033u_1 - 0.005u_5)(1.795u_2 - 2.109u_5) - (0.040u_1 - 0.042)\exp(1.895u_1)$\\
            &$g_4$ & \small $-0.999(0.009u_1 - 0.004u_5 + 0.009)(1.449u_2 + 3.357u_3 + 1.044u_5 + 2.355)$\\
            &$g_5$ & \cellcolor{lightgray} \small $-0.017u_1 - 0.013u_3^2 + 0.011u_5 + 0.852(0.001u_1 - 0.001u_5)(17.071u_4 - 31.025u_5) - 0.009$\\
			\bottomrule
		\end{tabular}
	\end{center}
 \label{tab:2dcyl_symbolic_expressions}
\end{table}
In Table~\ref{tab:2dcyl_symbolic_expressions}, we report the best symbolic models in terms of the $\text{rMSE}_{\text{test}}$ for both Reynolds numbers.
We observe that all the expressions are simple and compact.

\begin{figure}[!ht]
    \centering
    \includegraphics[scale=1]{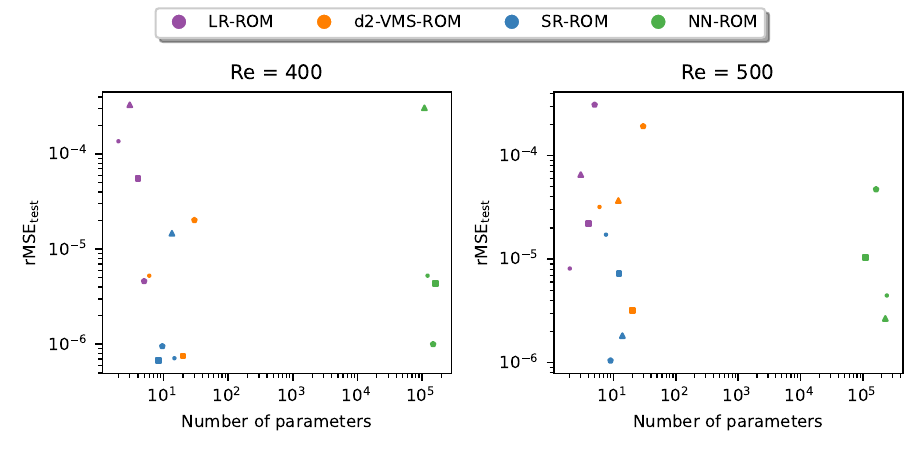}
    \caption{2D flow past a cylinder. $\text{rMSE}_{\text{test}}$ as a function of the number of parameters of the model. For the SR-ROM and the NN-ROM, the mean values out of 5 independent runs are reported. Each point is marked using the following legend: point for $r=2$; triangle for $r=3$; square for $r=4$; pentagon for $r=5$. 
    }
    \label{fig:pareto_2dcyl}
\end{figure}
In Figure~\ref{fig:pareto_2dcyl}, we further plot 
the $\text{rMSE}_{\text{test}}$ as a function of the number of parameters for each ROM. We note that the number of constants is reported for the SR-ROM and the number of weights for the other ROMs.
We remark that the d2-VMS-ROM results for $\rm{Re}=400$ and $r=3$ 
do not appear in the plot, since the model's 
kinetic energy diverges. Also, 
although the $\text{rMSE}_{\text{test}}$s of SR-ROM and NN-ROM  are similar for $\rm{Re}=500$  (Table~\ref{tab:2dcyl_test_error}), SR-ROM has an order of magnitude fewer parameters.
Overall, the SR-ROM is the best model in terms of complexity and accuracy.

\begin{figure}[!ht]
    \centering
    \includegraphics[width=\columnwidth]{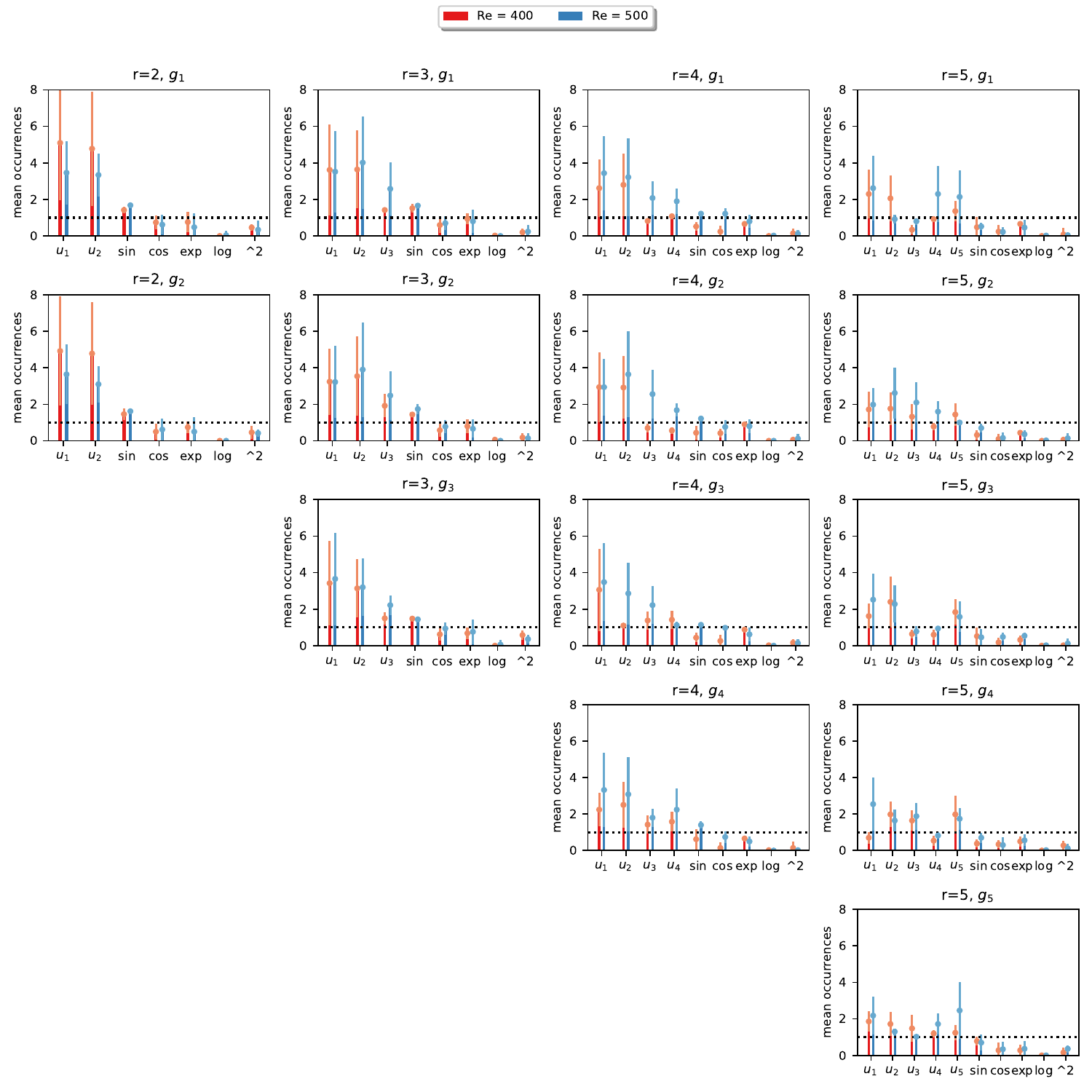}
    \caption{2D flow past a cylinder. 
    Statistics of the occurrences (over 50 independent runs) of the ROM coefficients 
    and primitives appearing in the function $g_i$ found using symbolic regression.
    }
    \label{fig:bar_plot_2dcyl}
\end{figure}

Finally, the structure of the 
closure components $g_i$ 
obtained using symbolic regression for both Reynolds numbers is 
analyzed in Figure~\ref{fig:bar_plot_2dcyl}.  
Specifically, we count the mean occurrences of the ROM coefficients, and square, trigonometric, and exponential functions over 50 independent runs after simplifying the symbolic expressions. 
To clarify 
our definition of ``occurrences'', let us consider the expression of the function $g_1$
for $\rm Re=500$, as shown in Table~\ref{tab:2dcyl_symbolic_expressions},
\begin{equation}
    g_1(u_1, u_2, u_3, u_4, u_5) := 0.015u_1 + 0.01u_4 + 1.0\sin(0.191u_1 - 0.1u_5 + \sin(0.064u_4) - 1.57) + 0.98.
\end{equation}
In the function $g_1$, there are two $u_1$, two $u_4$, one $u_5$, and two sine occurrences. 
The same rule 
is applied to 
all the components for each run, and we compute the 
{mean and standard deviation (std)} 
with respect to the number of runs.
From Fig.~\ref{fig:bar_plot_2dcyl}, we 
can observe that:
\begin{itemize}
    \item The logarithm and square are absent from almost all the expressions, 
    thus suggesting that the (general) primitive set that we 
    chose for SR could be reduced without loss in accuracy of the models;
    \item Except for $r=5$, the first ROM coefficients $u_1$ and $u_2$ appear more frequently than the others (even though their frequency has a significant standard deviation). 
    Intuitively, this could be related to the higher importance of the first modes in the reconstruction of the solution;
    \item For $r=5$, the mean occurrences of all the considered terms are much more spread than the previous $r$ 
    values.
\end{itemize}
We note that the expressions of $g_i$ found via SR do not show any common structure across the two values of Reynolds number. This is possibly a consequence of the training procedure, which was performed by 
including only data corresponding to a single Reynolds number.

\clearpage
\subsection{2D Lid-Driven Cavity}
\label{subsec: 2dldc}

Our next example is the 2D lid-driven cavity (LDC) problem, which is a more challenging test problem compared to the 2D flow past a cylinder. As demonstrated in \cite{fick2018stabilized}, the problem requires more than $60$ POD modes for G-ROM to accurately reconstruct the solutions and quantities of interest. 
A detailed description of the FOM setup for this problem can be found in \cite{kaneko2020towards}. 
\begin{figure}[!ht]
    \centering
     \includegraphics[width=0.75\textwidth]{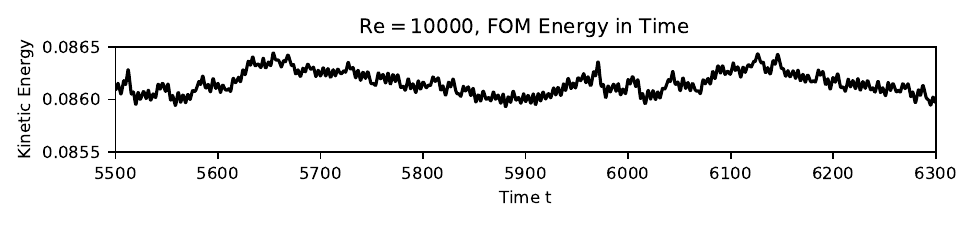}
     \includegraphics[width=0.75\textwidth]{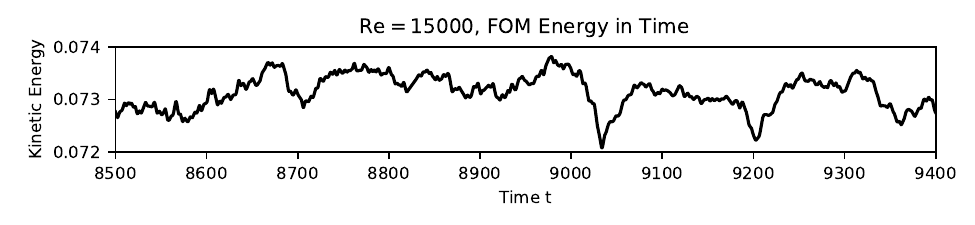}
     \includegraphics[width=0.75\textwidth]{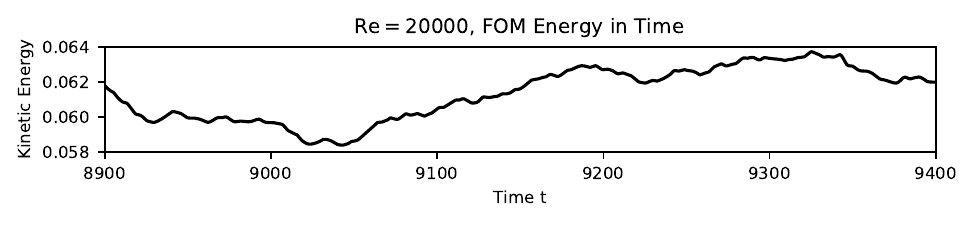}
    \caption{2D lid-driven cavity. The FOM kinetic energy {as a function of} time $t$ for $\rm Re=10000,15000,20000$.}
    \label{fig:2dldc_fom_energy}
\end{figure}
We consider three Reynolds numbers, namely, $\rm Re=10000$, $\rm Re=15000$, and
$\rm Re=20000$. 
The time evolution of the FOM kinetic energy for the three Reynolds numbers is shown in Fig.~\ref{fig:2dldc_fom_energy}. The results indicate that the solutions at these Reynolds numbers are chaotic.
\begin{figure}[!ht]
    \centering
     \includegraphics[width=0.95\textwidth]{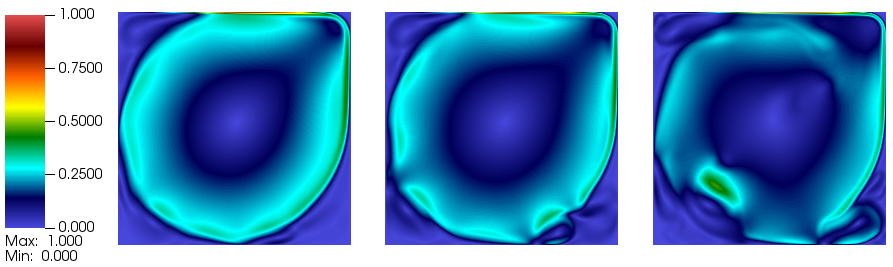}
    \caption{2D lid-driven cavity. The velocity magnitude snapshot for $\rm Re=10000,15000,20000$.} 
    \label{fig:2dldc_fom_snap}
\end{figure}

For each $\rm Re$, $K=2001$ snapshots $\{\bu^k :=
\bu(\bx,t^k)-\bphi_0\}^K_{k=1}$ are collected in the time interval
$[T_0,~T_0+80]$ with sampling time step $\Delta t_s=0.04$, and the zeroth mode,
$\bphi_0$, is set to be the FOM velocity field at $t=T_0$. We note that the time
$T_0$ depends on the 
Reynolds number and is determined such that the solutions are in
the statistically steady state region.
Snapshots of the velocity magnitude for each considered $\rm Re$ are shown in Fig.~\ref{fig:2dldc_fom_snap}.
The snapshot data are then used to construct the reduced basis functions
$\{\bphi_i\}^r_{i=1}$ by 
leveraging POD, the VMS-ROM closure matrix $\bTau$, and \eqref{eq:vmsrom-closure}. The data-driven closure of each model is trained in the time interval $[T_0,~T_0+80]$ by minimizing \eqref{eq: mse_tr}, and is selected based on the optimal hyperparameters that minimize~\eqref{eq: mse_val} in the time interval $[T_0+80,~T_0+160]$.
Finally, we assess the generalization capabilities of all the ROMs in the extrapolation region, that is, in the time interval $[T_0,~T_0+400]$.

\begin{table}[!ht]
    \centering
    \caption{2D lid-driven cavity. Relative MSE ($\text{rMSE}_{\text{test}}$) between the FOM kinetic energy and the predicted kinetic energy in the testing region for all models and reduced space dimension $r=2,3,4,5,6,7$. For the SR-ROM and the NN-ROM, the mean and the standard deviation over $5$ independent runs (different seeds) are reported. 
    In bold, the best performing model for each case is listed.
    }
    \begin{tabular}{ccccccc}
         \toprule
          Re & $r$ & LR-ROM & d2-VMS-ROM & SR-ROM & NN-ROM\\
          \midrule
          \multirow{6}{*}{10000} & 2& $2.8 \times 10^{-5}$ & $2.8 \times 10^{-5}$  & $(2.8 \pm 0.0) \times 10^{-5}$ & \boldmath $(4.1 \pm 4.2) \times 10^{-6}$ \\
          & 3& $0.2$ & $2.8 \times 10^{-3}$  & $(2.2 \pm 2.2) \times 10^{-5}$ & \boldmath $(2.9 \pm 1.3) \times 10^{-6}$\\
          & 4 & $2.2$ & $0.7$  & \boldmath $(4.9 \pm 1) \times 10^{-6}$ & \boldmath $(5.0\pm 3.3) \times 10^{-6}$\\
          & 5 & $1.5 \times 10^{-3}$ & $1.5 \times 10^{-3}$  & \boldmath $(4.2 \pm 5.2) \times 10^{-5}$ & $(5.4 \pm 1.5) \times 10^{-3}$\\
          & 6 & $3.1 \times 10^2$ & $7$  & \boldmath $(1.5 \pm 1.3) \times 10^{-5}$ & $(2.7 \pm 5.5) \times 10^2$\\
          & 7 &\boldmath $2.3 \times 10^{-5}$ & \boldmath $2.5 \times 10^{-5}$ & \boldmath $(2.9 \pm 0.2) \times 10^{-5}$ &\boldmath $(1.9 \pm 3.1) \times 10^{-5}$\\
          \midrule
          \multirow{6}{*}{15000} & 2 & $1.4 \times 10^{-4}$ & \boldmath $8.4 \times 10^{-5}$ & $(2.8 \pm 0.6) \times 10^{-4}$ & \boldmath $(4.5 \pm 4.3) \times 10^{-5}$\\
          & 3 & $5.1 \times 10^{-3}$ & $1.8 \times 10^{-4}$ & \boldmath $(3.9 \pm 0.5) \times 10^{-5}$ & $(0.6 \pm 1.2) \times 10^6$\\
          & 4 & $4 \times 10^{-4}$ & $1.6 \times 10^{-4}$ & \boldmath $(8.8 \pm 1.5) \times 10^{-5}$ & \boldmath $(1.0 \pm 0.4) \times 10^{-4}$\\
          & 5 & $1.9 \times 10^{-4}$ & \boldmath $9.2 \times 10^{-5}$ & $0.01 \pm 0.01$ & \boldmath $(7.2 \pm 3.8) \times 10^{-5}$\\
          & 6 & $0.06$ & $1.4 \times 10^{-4}$ & \boldmath $(4.8 \pm 0.8)\times 10^{-5}$ & $(6.6 \pm 10) \times 10^{-4}$\\
          & 7 & $0.12$ & $1.7 \times 10^{-4}$ & \boldmath $(8.2 \pm 3.3) \times 10^{-5}$ & $0.04 \pm 0.07$\\
          \midrule
          \multirow{6}{*}{20000} & 2 & $1.5 \times 10^{-3}$ &  $1.1 \times 10^{-3}$ & \boldmath $(2.2 \pm 0.6) \times 10^{-4}$ & $(1.7 \pm 0.1) \times 10^{-3}$\\
          & 3 & \boldmath $1.1 \times 10^{-4}$ & $6.3 \times 10^{-4}$ & $(1.7 \pm 0.4) \times 10^{-4}$ & $(5.1 \pm 6.3) \times 10^{-4}$\\
          & 4 & 0.3 & $6.2 \times 10^{-4}$ & \boldmath $(1.8 \pm 0.4) \times 10^{-4}$ & $0.7 \pm 0.7$\\
          & 5 & $2.6 \times 10^{-2}$ &  \boldmath $3.6 \times 10^{-3}$ & \boldmath $(3.3 \pm 6.2) \times 10^{-3}$ & $(3.8 \pm 5.2) \times 10^{-2}$\\
          & 6 & \boldmath $8.1 \times 10^{-4}$ & $2.9 \times 10^{-3}$ & nan & \boldmath $(1.8 \pm 1.3) \times 10^{-3}$\\
          & 7 & \boldmath $2 \times 10^{-4}$ & $8.6 \times 10^{-4}$ & \boldmath$(3.4 \pm 2.6) \times 10^{-4}$ & $0.1 \pm 0.2$\\
          \bottomrule
    \end{tabular}
    \label{tab:ldc_test_error}
\end{table}
In Table~\ref{tab:ldc_test_error}, the $\text{rMSE}_{\text{test}}$ between the FOM kinetic energy and the predicted kinetic energy in the extrapolation region of all models 
is reported with respect to reduced space dimension $r=2,3,4,5,6,7$. The rMSE of the SR-ROM and the NN-ROM is computed by taking the mean of $5$ runs with $5$ different seeds because these models rely on stochastic optimization algorithms. 
Similar to the 2D flow past a cylinder case, the results show that the SR-ROM is the most robust model, followed by the NN-ROM, d2-VMS-ROM, and LR-ROM. {In particular, we note that for $\rm{Re} = 10000$ and $r=5,6$, $\rm{Re} = 15000$ and $r=3,6,7$, and $\rm{Re} = 20000$ and  $r=2$, the SR-ROM $\text{rMSE}_{\text{test}}$ in the extrapolation region is orders of magnitude lower than the $\text{rMSE}_{\text{test}}$ of the other ROMs.}

We further measure the closeness between the VMS-ROM closure matrix $\bTau$ and the reconstructed closure for each model {by computing the mean MSE~\eqref{eq: mean_mse_tr}, see Table~ \ref{tab:ldc_tau_rec_error}.} 
We observe that the SR-ROM is the most robust method also in this case. Furthermore, we note that, unlike in the MSE of predicted energy (Table~\ref{tab:ldc_test_error}), the NN-ROM has a much larger MSE than the SR-ROM. {We stress again that this phenomenon is a consequence of the $L^2$ regularization, which is not used for the SR-ROM.}
\begin{table}[!ht]
    \centering
    \caption{2D lid-driven cavity. Mean squared error ({$\overline{\text{MSE}_{\text{tr}}}$}) between the VMS-ROM closure and the ML-VMS-ROM reconstructed closures for $\rm Re=10000,~15000,~20000$. In bold, the best performing model for each case is listed.}
    \begin{tabular}{ccccccc}
         \toprule
          Re & $r$ & LR-ROM & d2-VMS-ROM & SR-ROM & NN-ROM\\
          \midrule
          \multirow{6}{*}{10000} & 2 & \boldmath $2.0 \times 10^{-7}$ & \boldmath $2.0 \times 10^{-7}$ & \boldmath $(1.0 \pm 0.0) \times 10^{-7}$& $(8.9 \pm 3.7) \times 10^{-4}$\\
          & 3 & \boldmath $1.4 \times 10^{-7}$ & \boldmath $1.2 \times 10^{-7}$ & \boldmath $(9.4 \pm 0.5) \times 10^{-8}$ & $(3.3 \pm 2.1) \times 10^{-4}$\\
          & 4 & $1.2 \times 10^{-7}$ & $1.4 \times 10^{-7}$ & \boldmath $(7.2 \pm 0.8) \times 10^{-8}$ & $(4.7 \pm 2.6) \times 10^{-4}$\\
          & 5 & $1.3 \times 10^{-7}$ & $1.3 \times 10^{-7}$ & \boldmath $(6.4 \pm 0.8) \times 10^{-8}$ & $(2.4 \pm 0.9) \times 10^{-7}$\\
          & 6 & $9.4 \times 10^{-8}$ & $9.3 \times 10^{-8}$ & \boldmath $(5.8 \pm 0.6) \times 10^{-8}$ & $(1.9 \pm 1.3) \times 10^{-4}$\\
          & 7 & $1.0 \times 10^{-7}$ & $1.0 \times 10^{-7}$ & \boldmath $(5.2 \pm 0.7) \times 10^{-8}$ & $(1.2 \pm 1.2) \times 10^{-5}$\\
          \midrule
          \multirow{6}{*}{15000} & 2 & $1.3 \times 10^{-5}$ & $1.4 \times 10^{-5}$ & \boldmath $(6.6 \pm 1.8) \times 10^{-6}$ & $(4.8 \pm 1.1) \times 10^{-4}$\\
          & 3  & $7.9 \times 10^{-6}$ & $4.5 \times 10^{-6}$ & \boldmath $(3.5 \pm 0.0) \times 10^{-6}$ & $(1.9 \pm 2.5) \times 10^{-4}$\\
          & 4 & $5.3 \times 10^{-6}$ & $3.1 \times 10^{-6}$ & \boldmath $(2.9 \pm 0.4) \times 10^{-7}$ & $(0.1 \pm 2.2) \times 10^{-3}$\\
          & 5 & $4.0 \times 10^{-6}$ & \boldmath $2.3 \times 10^{-6}$ & \boldmath $(2.3 \pm 0.4) \times 10^{-6}$ & $(9.7 \pm 7.1) \times 10^{-6}$\\
          & 6 & $3.6 \times 10^{-6}$ & \boldmath $2.1 \times 10^{-6}$ & \boldmath $(2 \pm 0.4) \times 10^{-6}$ & $(8.7 \pm 3.8) \times 10^{-6}$\\
          & 7 & \boldmath $2.7 \times 10^{-6}$ & \boldmath $2.5 \times 10^{-6}$ & \boldmath $(1.7 \pm 1.9) \times 10^{-6}$ & $(1.1 \pm 0.5) \times 10^{-5}$\\
          \midrule
          \multirow{6}{*}{20000} & 2 & $9.0 \times 10^{-6}$ & $9 \times 10^{-6}$ & \boldmath  $(5.4 \pm 0.6) \times 10^{-6}$ & $(9.0 \pm 0.0) \times 10^{-6}$\\
          & 3 & $9.2 \times 10^{-6}$ & $8.7 \times 10^{-6}$ & \boldmath  $(4.7 \pm 1.2) \times 10^{-6}$ & $(4.2 \pm 5.5) \times 10^{-6}$\\
          & 4 & $1.1 \times 10^{-5}$ & $9.3 \times 10^{-6}$ & \boldmath  $(4.3 \pm 0.5) \times 10^{-6}$ & $(1.1 \pm 0.1) \times 10^{-5}$\\
          & 5 & $9.8 \times 10^{-6}$ & $1.0 \times 10^{-5}$ & \boldmath  $(6.5 \pm 0.9) \times 10^{-6}$ & $(3.2 \pm 4) \times 10^{-5}$\\
          & 6 & $7.4 \times 10^{-6}$ & $7.6 \times 10^{-6}$ & \boldmath  $(4.9 \pm 0.5) \times 10^{-6}$ & $(1.0 \pm 1.8) \times 10^{-4}$\\
          & 7 &\boldmath  $7.3 \times 10^{-6}$ & \boldmath  $7.6 \times 10^{-6}$ &\boldmath  $(5.3 \pm 0.8) \times 10^{-6}$ & $ \boldmath (7.6 \pm 0.1) \times 10^{-6}$\\
          \bottomrule
    \end{tabular}
    \label{tab:ldc_tau_rec_error}
\end{table}

\begin{figure}[!ht]
    \centering
    \includegraphics[scale=0.52]{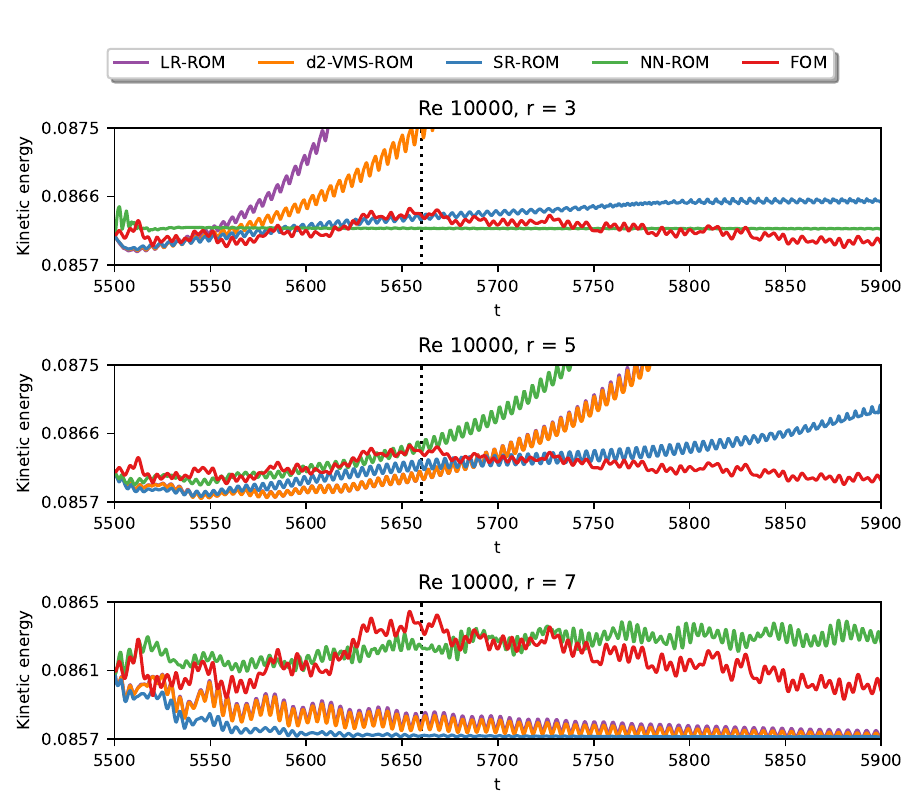}
     \includegraphics[scale=0.52]{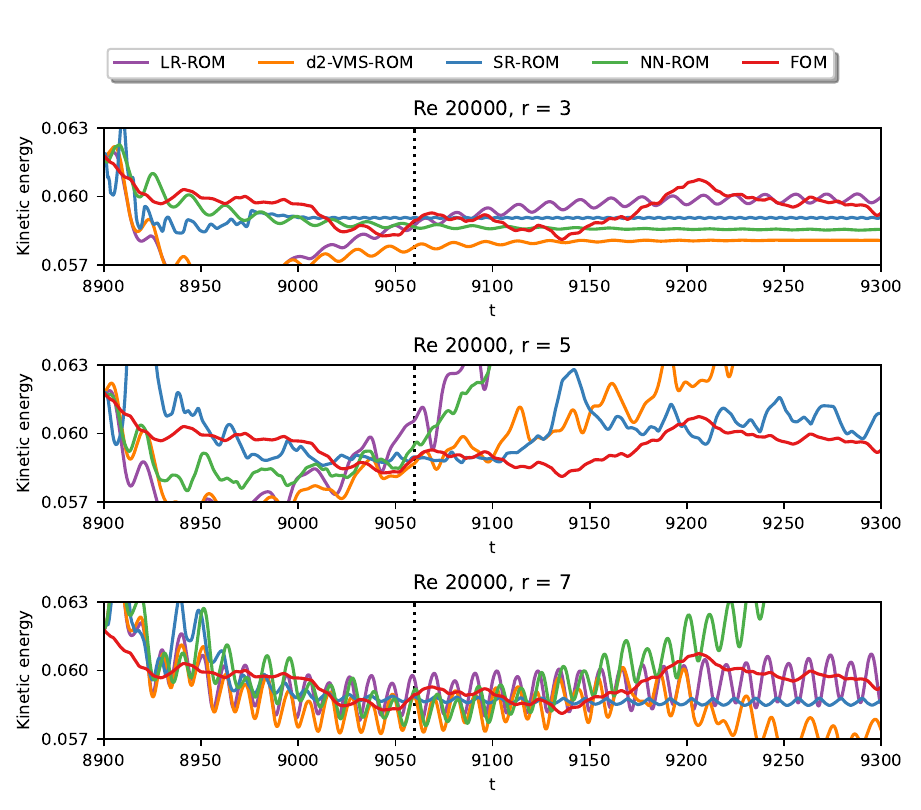}
    \caption{2D lid-driven cavity at $\rm Re=10000$ and $\rm Re=20000$.  Comparison of the kinetic energy between the FOM and 
    the ROMs with $r=3,5,7$.}
    \label{fig:energies_ldc}
\end{figure}

Similarly to Section~\ref{subsec: 2dfpc}, we report in Table~\ref{tab:ldc_symbolic_expressions} the best symbolic expressions found for the three Reynolds numbers considered. 
These expressions are slightly more complicated than the ones in Table~\ref{tab:2dcyl_symbolic_expressions}, but this is due to the chaotic, thus more complex, nature of the problem. Also, we remark that for $\rm Re= 15000$ and $20000$, the best models found depend only polynomially on the FOM coefficients $u_r$. Hence, also the d2-VMS-ROM could, in principle, find the same models. However, with respect to the $\text{rMSE}_{\text{test}}$ in the extrapolation region, a worse model is found, and this is due to the ill-conditioned linear system to be solved to train this method \cite{xie2018data}. 

\begin{figure}[!t]
    \centering
    \includegraphics[scale=1.05]{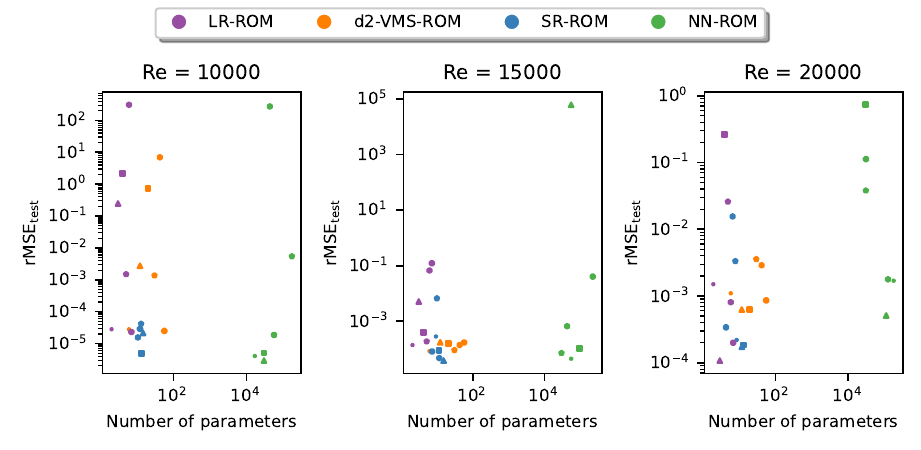}
    \caption{2D lid-driven cavity. $\text{rMSE}_{\text{test}}$ as a function of the number of parameters of the model. For the SR-ROM and NN-ROM, the mean values over 5 independent runs are reported. Each point is marked using the following legend: point for $r=2$; triangle for $r=3$; square for $r=4$; pentagon for $r=5$; hexagon for $r=6$; octagon for $r=7$.}
    \label{fig:pareto_ldc}
\end{figure}

Figure~\ref{fig:pareto_ldc} confirms, also in this case, that the SR-ROM is 
the best trade-off between generalization capability ($\text{rMSE}_{\text{test}}$) and model complexity. Again, the SR-ROM and NN-ROM manifest a similar behavior in terms of $\text{rMSE}_{\text{test}}$, especially for $\rm{Re}=10000$ and $15000$, but the SR-ROM has an order of magnitude fewer parameters. Also, the SR-ROM results for $\rm{Re} = 20000$ and $r=6$ are not shown in the figure, since the kinetic energy of one of the SR-ROM best models diverges.

As in Section~\ref{subsec: 2dfpc}, we compute the mean occurrences of some terms using the same strategy. Also in this case, we can make a qualitative analysis:
\begin{itemize}
    \item The logarithm and square are basically absent for every $r$ and $\rm Re$;
    \item For $r=5,6,7$ and for each Reynolds number, the mean occurrence of every term is uniform, \textsl{i.e.},~none of the terms in the expressions is dominant. 
\end{itemize}
As in Section~\ref{subsec: 2dfpc}, 
no evident structure emerges in the expressions of  Table~\ref{tab:ldc_symbolic_expressions} across different values of the Reynolds number because the training procedure of all the ROMs is 
carried out at a fixed Reynolds number.

\begin{table}[!t]
    \caption{2D lid driven cavity. 
    Expressions of the components of the SR-ROM closure model that attains the minimum $\text{rMSE}_{\text{test}}$ in the extrapolation region 
    among those reported in Table~\ref{tab:ldc_test_error}.} 
	\begin{center}
		\begin{tabular}{ccc}
			\toprule
			Re & Component &  Expression \\
			\midrule
            \multirow{11}{*}{10000} & \multirow{3}{*}{$g_1$} & \cellcolor{lightgray} \small $0.0002u_1 + 0.0901u_2 + 0.0003u_3 - 0.0004u_4 +(0.0001u_1 - 0.0001u_2)* $\\
            & & \cellcolor{lightgray} \small $*\cos(3.7331u_4 - 2.1130) + (0.0002u_1 + 0.0010u_2 + 0.0005u_3 - 0.0007u_4 + 0.0004)*$\\
            & & \cellcolor{lightgray} \small $*\sin(0.0686u_1 - \cos(0.6297u_1 - 0.1635u_3)) - 0.999\sin(0.0897u_2) + 0.0003$\\
            &\multirow{3}{*}{$g_2$} & \small $-0.0003u_1 + 0.0003u_2 + 0.0002u_3 - \sin((0.0883u_2(1.7855u_1 + 1.8087u_2) - $\\
            & & \small $- \sin(1.2270u_2 + 0.5128u_4 - 7.0881) + 0.9256)\sin(\sin(\sin(0.0006u_3) - 0.0009))*$\\
            & & \small $*\sin(\sin(\sin(\sin(0.8346u_1)) + \cos(0.7177u_1 + 0.9317)))) - 0.0004$\\
            &\multirow{3}{*}{$g_3$} & \cellcolor{lightgray} \small $0.0001u_2 + 0.0001u_3 + (0.0005u_1u_3 - 0.0003u_2 - 0.0003u_4)(0.4274u_1 - 0.1150u_4 + $\\
            & &\cellcolor{lightgray} \small $+ (0.0507u_1 - 0.0376u_2\cos(1.1193u_2) + 0.0342u_3)(1.6577u_1 + 0.5313u_2 - 7.1728) + $\\
            & &\cellcolor{lightgray} $ + \sin(0.7269u_2 + 1.0911u_3 + 0.5615)\cos(0.6477u_1))$ \\
            &\multirow{3}{*}{$g_4$} & \small $1.0010\sin(0.0299u_2 + 0.0659u_4 + (0.3499u_1 - \cos(1.6868u_1))(0.4145u_4 + \sin((0.4721u_2 -$\\
            & & \small $- \sin(0.7269u_3))\cos(0.6962u_1) + 2.7152)0.0006 -\sin(0.0301u_2 + 0.0659u_4))$\\
            \hline
            \multirow{8}{*}{15000} &\multirow{2}{*}{$g_1$} &\cellcolor{lightgray} \small $-0.0022 - 0.0006u_1 + 0.0004u_0^2(1.9712 + 0.5332u_0) - 0.0001u_1u_0 - $\\
            & &\cellcolor{lightgray} \small $- (-0.0004 - 0.0001u_1)(0.7896u_2 + 2.7564u_1^2)$\\
            &\multirow{3}{*}{$g_2$} & \small $0.0001 + 0.9990(20.9826 + 9.6691u_0 + 7.9730u_2 + (1.1479 + 0.0092u_0(-22.8110u_0 $\\
            & & \small $+ 27.8969u_1 + 9.8753u_2) - 3.1177u_1(0.1684 + 0.0210u_2))(-13.6013 - 8.9292u_1 + $\\
            & & \small $+(-2.5189 + 0.2241u_0 - 1.2826u_2)(6.0625 +  2.3368u_0 + 1.1539u_1 - 0.0070u_2)))$\\
            & & \small $(-0.0002 + 0.0001u_1 + 0.0002u_2)$\\
            &\multirow{2}{*}{$g_3$} &\cellcolor{lightgray} \small $0.0038 - 0.0006u_0 - 0.0019u_2 + 0.0001u_2u_0 + (0.0010 - 0.0007u_2)(-2.4952  +$\\
            & &\cellcolor{lightgray} \small $+ 0.0602u_1 + 0.1949u_2 + (0.2424u_1 + 0.1948u_2u_1)(1.1255u_0 - 1.2423u_1 + 2.9477u_2))$\\
            \hline
            \multirow{10}{*}{20000} & \multirow{2}{*}{$g_1$} & \small $0.0002u_1 + 0.0008u_2 - 0.0005u_3 + (6.0344u_2 - 4.848u_4)(0.0001u_1(1.9919u_1u_3*$\\
            & & \small $*(3.2164u_1 + 1.9255u_4 - 3.5296) - 4.4491u_2 - 2.9491u_3 + 8.9232)) - 0.0004$\\
            &\multirow{3}{*}{$g_2$} &\cellcolor{lightgray} \small $1.0021(0.3193u_2(0.7430u_1 + 1.8020u_4)(0.1953u_2 - 1.4475)(5.9033u_2 + 10.0442) + $ \\
            & &\cellcolor{lightgray} \small $+ 12.5300u_4 - (0.5016u_2 - 4.5800u_3)(1.1975u_2 + 2.6217)(0.3743u_3 - 0.0664))*$\\ 
            & &\cellcolor{lightgray} \small $*(0.0002u_1 + 0.0001u_2) - 0.0001$\\
            &\multirow{2}{*}{$g_3$} &  \small $0.0003u_3 - 0.0002u_4 + (0.0001u_2 - 0.0001u_3)(1.7047u_2^2 + 2.8970u_4 - 0.9176)(13.6070u_1 - $\\
            & &  \small $ - 1.5597u_3 + 5.6265u_4 - 6.6889) - (0.0021u_2 + 0.0011u_4)(1.5428u_2 - 0.6498u_3) + 0.002$\\
            &\multirow{3}{*}{$g_4$} &\cellcolor{lightgray} \small $-0.9953(0.0001u_2 + 0.0001u_3^2 - 0.0001u_3 + 0.001u_4 - (0.0001u_1 + 0.0010)*$\\
            & &\cellcolor{lightgray} \small $*(1.2730u_4 + 0.300))(4.1908u_2^2 + 9.7453u_2 - 0.1684u_3(70.4595u_1 +$\\
            & &\cellcolor{lightgray} \small $+  13.7008) + 4.0931u_4 + 3.2191)$\\
			\bottomrule
		\end{tabular}
	\end{center}
 \label{tab:ldc_symbolic_expressions}
\end{table}

In Figure~\ref{fig:energies_ldc}, we show the kinetic energy 
as a function of time for the ROMs and the FOM with reduced space dimension $r=3,5,7$, for $\rm Re=10000$ and $\rm Re=20000$. 
For $\rm Re=10000$, the SR-ROM outperforms the other ROMs for $r=5$, while for $r=3,7$ the SR-ROM is still among the best 
performing models. We observe that the d2-VMS-ROM and LR-ROM deviate in the extrapolation region for $r=3,5$, as does the NN-ROM for $r=5$. The SR-ROM is the only method that 
remains close to the FOM in this case. Finally, for $\rm Re=20000$, the SR-ROM is the best method for $r=3,5,7$.

\begin{figure}[!t]
    \centering
    \includegraphics[width=\columnwidth]{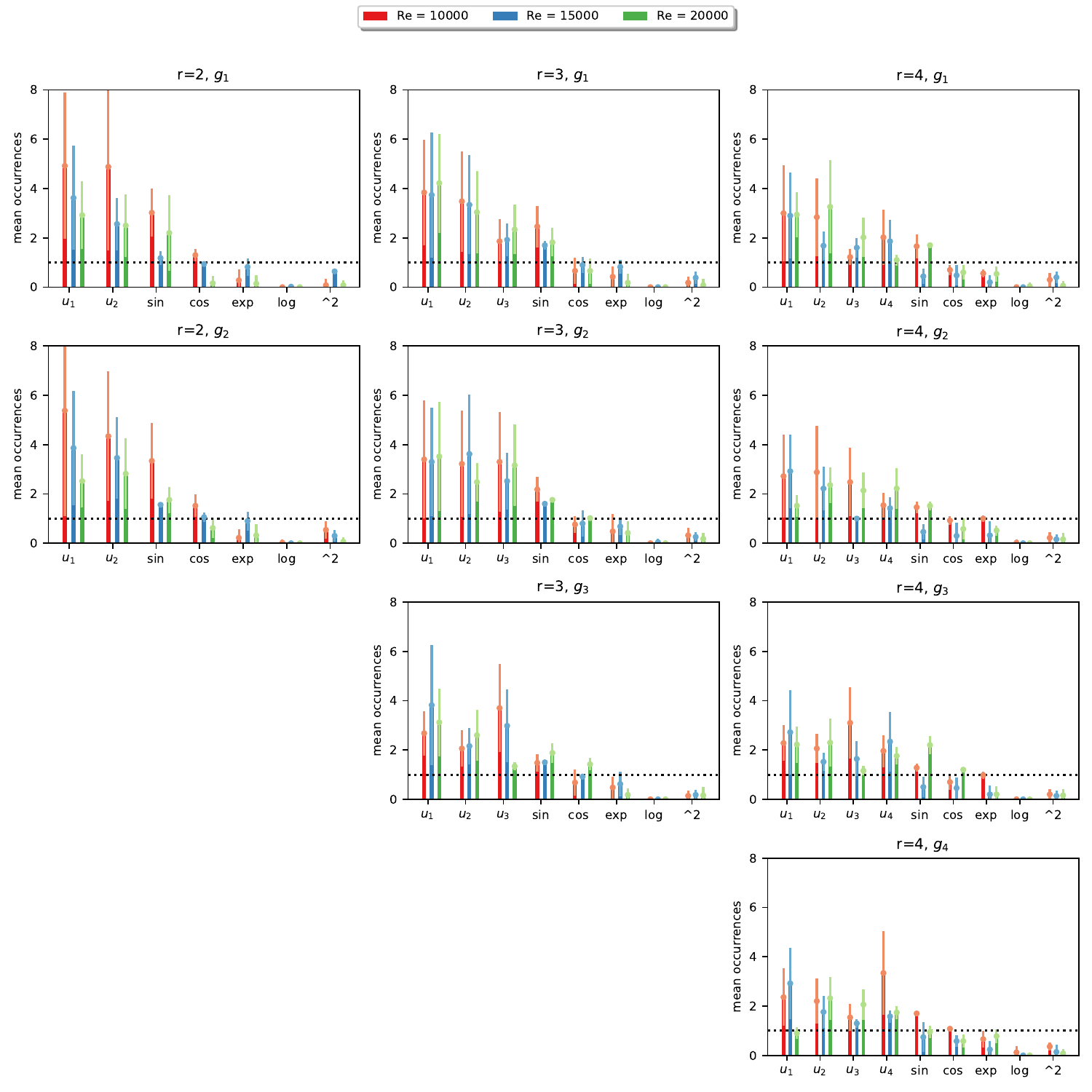}
    \caption{2D lid-driven cavity, $r=2,3,4$. 
    Statistics of the occurrences (over 50 independent runs) of the ROM coefficients 
    and primitives appearing in the closure model found by SR-ROM.}
    \label{fig:bar_plot_ldc_2_3_4}
\end{figure}

\begin{figure}[!t]
    \centering
    \includegraphics[width=\columnwidth]{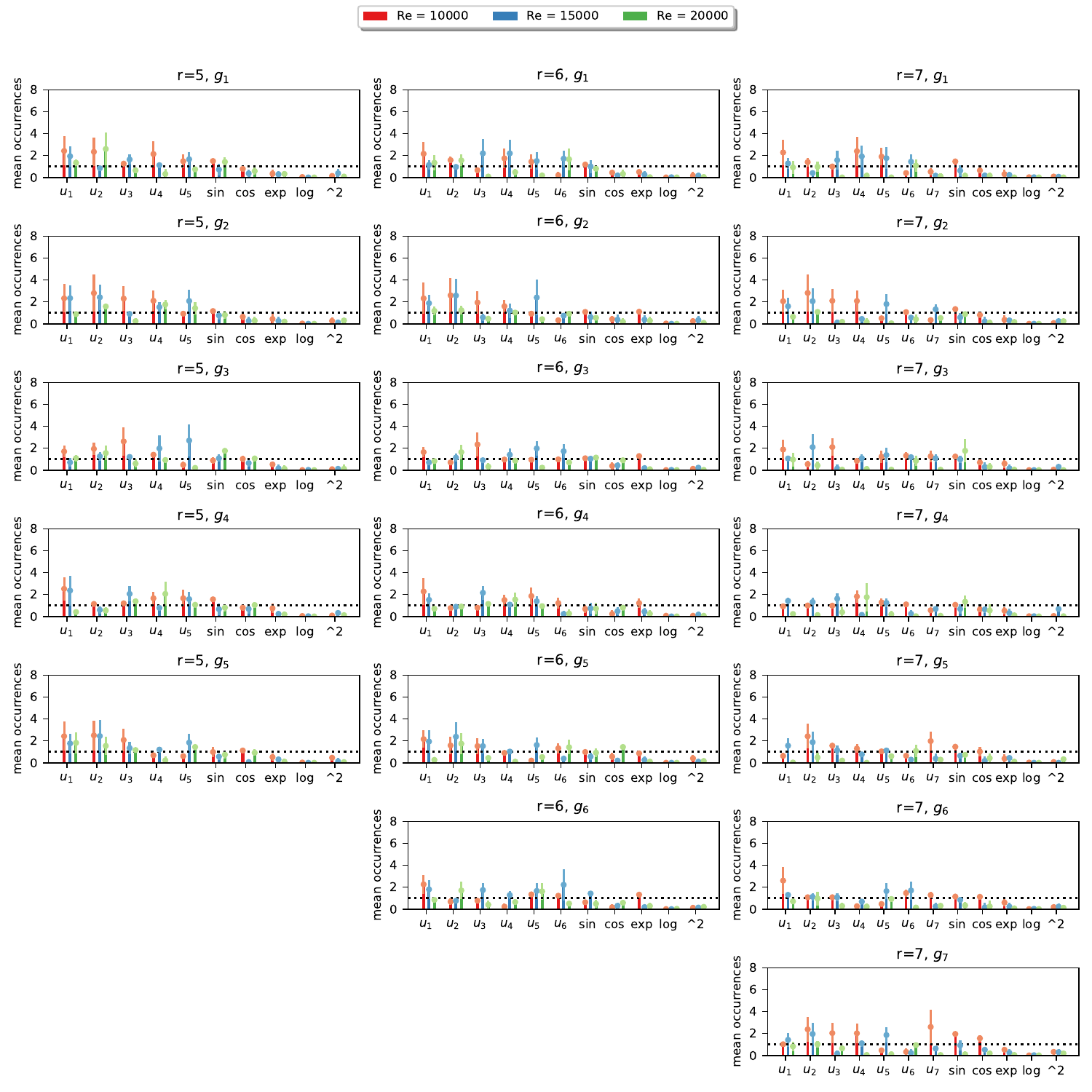}
    \caption{2D lid-driven cavity, $r=5,6,7$. 
    Statistics of the occurrences (over 50 independent runs) of the ROM coefficients 
    and primitives in the closure model found by SR-ROM.}
    \label{fig:bar_plot_ldc_5_6_7}
\end{figure}

\clearpage

\section{Conclusions and Future Work}
\label{sec:conclusions}

In this paper, we propose a novel symbolic regression (SR) closure strategy for ROM, which leverages available data and SR to find the ROM closure model.
The new SR-ROM closure model combines the advantages of current data-driven ROM closure approaches and eliminates their drawbacks:
The new data-driven SR closures yield ROMs that are interpretable, parsimonious, accurate, generalizable, and robust.
To compare the SR-ROM closures with the current data-driven ROM closures (i.e., structural and machine learning ROM closures), we consider a specific ROM closure framework, the data-driven variational multiscale ROM.
Within this framework, we compare the novel SR-ROM closure with three other data-driven ROM closures: 
a linear regression ROM closure, a quadratic data-driven ROM closure, and a neural network ROM closure.
To investigate the four types of ROM closures, we consider two under-resolved, convection-dominated test problems: the flow past a cylinder and the lid-driven cavity flow.
To assess the performance of the different types of closures, we monitor the performance of the  ROMs 
obtained by adding the four ROM closure models. 

Our numerical investigation yields the following conclusions: 
\begin{itemize}
\item For both test problems, in the predictive regime, the new SR-ROM closure yields more accurate results than the other three ROM closures (i.e., the linear regression, quadratic data-driven, and neural network ROM closures);
\item The SR-ROM is the most robust model with respect to the ROM space dimension, $r$;
\item Our numerical results also show that the new SR-ROM closure has orders of magnitude fewer terms than the neural network closure model;
\item The SR strategy is more accurate than the other three data-driven approaches in reconstructing the ROM closure term.  
\end{itemize}

The first steps in the numerical investigation of the new SR-ROM have been encouraging. There are, however, several other research directions that should be pursued next. For example, in this paper, we 
did not find common structures 
in the SR-ROM closures 
for different Reynolds numbers. This is because we 
used training data 
for a fixed Reynolds number. 
To allow SR to discover common structures in the new SR-ROM closures, 
we plan to train it with data from different Reynolds numbers. 
We also note that, although the SR-ROM is used in this paper for prediction in time, the Reynolds number is fixed. We 
plan to extend the new SR-ROM to 
parametric problems, 
and use it to predict solutions and quantities of interests 
at other Reynolds numbers.
Another interesting research direction is the investigation of the new SR-ROM strategy in ROM closure frameworks that are different from the variational multiscale ROM closure approach considered in this paper.
For example, we plan to investigate the new SR-ROM closure for large eddy simulation ROMs constructed by using ROM spatial filtering \cite{xie2018data}.
Finally, we plan to compare the new SR-ROM closure strategy with other data-driven modeling approaches, e.g., operator inference \cite{peherstorfer2016data,aretz2024exploiting}, SINDy \cite{brunton2016discovering,messenger2021weak}, and LaSDI \cite{fries2022lasdi}.

\section{CRediT Authorship Contribution Statement}

\textbf{Simone Manti:}
Writing - original draft, Writing - review and editing,
Data curation, Investigation, Conceptualization, Validation, Visualization, Methodology, Software.
\textbf{Ping-Hsuan Tsai:}
Writing - original draft, Writing - review and editing,
Data curation, Investigation, Conceptualization, Validation, Visualization, Methodology, Software.
\textbf{Alessandro Lucantonio:}
Writing - review and editing, 
Conceptualization, Validation, Methodology, Resources, Funding acquisition, Supervision.
\textbf{Traian Iliescu:}
Writing - original draft,
Writing - review and editing, Conceptualization, Validation, Methodology, Funding acquisition, Supervision. 

\section*{Acknowledgments}

{The work of SM and AL is supported by the European Union (European Research Council (ERC), ALPS, 101039481). Views and opinions expressed are however those of the author(s) only and do not necessarily reflect those of the European Union or the ERC Executive Agency. Neither the European Union nor the granting authority can be held responsible for them. Computational resources have been partially provided by DeiC National HPC (DeiC-AU-N1-2023030).}
{TI acknowledges support through National Science Foundation grants DMS-2012253 and CDS\&E-MSS-1953113.
}

\bigskip 

\bibliographystyle{unsrt}

\end{document}

%% file: pf.tex
\def\bX{{\bf X}}

\def\utau{{\bf{\underline \tau}}}

\def\uhu{{\widehat {\underline u}}}

\def\cN{{\cal N}}

\def\scriptO{{{\it O}\kern -.42em {\it `}\kern + .20em}}
\def\RR{{{\rm l}\kern - .15em {\rm R} }}
\def\PP{{{\rm l}\kern - .15em {\rm P} }}
\def\L2{{{\sf L}^2}}
\def\H1{{{\sf H}^1}}
\def\PN2{{\PP_{N}-\PP_{N-2}}}

\def\complex{{{\rm C} \kern - .53em {\rm l} \kern + .38em}}

\def\a1{{ | \lambda_{\min} |}}

\def\l1{{   \lambda_{\min}  }}

\def\bu0{{\underline {\bf 0}}}

\def\bu{{\bf u}}

\def\bv{{\bf v}}
\def\bx{{\bf x}}

\def\Oh{{\hat \Omega}}

\def\uf{{\underline f}}

\def\uu{{\underline u}}

\def\ux{{\underline x}}

\def\u0{{\underline 0}}
\def\1u{{\underline 1}}

\newcommand{\pp}[2]{\frac{\partial #1}{\partial #2} }

\def\tr{{\tilde r}}

\def\tR{{\! \tilde R}}

\def\cN{{\cal N}}

\def\mC{{\mathcal{C}}}

\def\mC{\mathcal{C}}

\def\bu{{\bm u}}
\def\bv{{\bm v}}

\def\tr{r}
\def\tR{R}

\newcommand{\bphi}{\boldsymbol{\varphi}}

\def\utau{\underline{\tau}}
\def\bTau{\bm{\mathcal{T}}}

\def\mP{{\mathcal{P}}}

\def\Ntr{N_{\text{tr}}}